\documentclass[fontsize=12pt,a4paper,headings=normal,
twoside=false,leqno,parskip=half-,abstract=true]{scrartcl}
\usepackage[english]{babel}
\usepackage[utf8]{inputenc}
\usepackage{hyperref}
\hypersetup{
 pdftitle={SturmForQuasilinearEq},
 pdfauthor={Phillipo Lappicy},
 colorlinks=true,
 linkcolor=blue,
 citecolor=blue,
 filecolor=blue,
 urlcolor=blue}

\usepackage{graphicx}
\usepackage[format=plain,labelfont=bf,font=small]{caption}
\usepackage{xcolor}
\usepackage[arrow, matrix, curve]{xy}
\usepackage{float}

\usepackage{caption}
\usepackage{tabulary}
\usepackage{array}
\newcolumntype{N}[1]{>{\centering\arraybackslash}m{#1}}

\usepackage{amsmath,amsthm}
\usepackage{amssymb} 

\makeatletter
\newcommand{\tpitchfork}{%
  \vbox{
    \baselineskip\z@skip
    \lineskip-.52ex
    \lineskiplimit\maxdimen
    \m@th
    \ialign{##\crcr\hidewidth\smash{$-$}\hidewidth\crcr$\pitchfork$\crcr}
  }%
}
\makeatother
\usepackage{latexsym}

\usepackage[notref,notcite,color,final 
]{showkeys}

\definecolor{refkey}{rgb}{1,0,0}
\definecolor{labelkey}{rgb}{1,0,0}

\usepackage{tikz}

\usepackage[textwidth=2cm,textsize=small,backgroundcolor=none]{todonotes}

  \mathchardef\ordinarycolon\mathcode`\:
  \mathcode`\:=\string"8000
  \begingroup \catcode`\:=\active
    \gdef:{\mathrel{\mathop\ordinarycolon}}
  \endgroup

\newtheorem{thm}{Theorem}[section]
\newtheorem{lem}[thm]{Lemma}
\newtheorem{prop}[thm]{Proposition}
\newtheorem{cor}[thm]{Corollary}
\newenvironment{pf}[1][Proof]{\begin{trivlist}
\item[\hskip \labelsep {\bfseries #1}]}{\end{trivlist}}

\usepackage{comment}
\usepackage{caption}

\begin{document}

\title{{\LARGE{Sturm attractors for quasilinear parabolic equations with singular coefficients}}}

\author{
 \\
{~}\\
Phillipo Lappicy*\\
\vspace{2cm}}

\date{ }
\maketitle
\thispagestyle{empty}

\vfill

$\ast$\\
Instituto de Ciências Matemáticas e de Computação\\
Universidade de S\~ao Paulo\\
Avenida trabalhador são-carlense 400\\
13566-590, São Carlos, SP, Brazil\\


\newpage
\pagestyle{plain}
\pagenumbering{arabic}
\setcounter{page}{1}

\begin{abstract}

The goal of this paper is to construct explicitly the global attractors of parabolic equations with singular diffusion coefficients on the boundary, as it was done without the singular term for the semilinear case by Brunovsk\'y and Fiedler (1986), generalized by Fiedler and Rocha (1996) and later for quasilinear equations by the author (2017). In particular, we construct heteroclinic connections between hyperbolic equilibria, stating necessary and sufficient conditions for heteroclinics to occur. Such conditions can be computed through a permutation of the equilibria. Lastly, an example is computed yielding the well known Chafee-Infante attractor.

\ 


\end{abstract}

\section{Main results}\label{sec:intro}

\numberwithin{equation}{section}
\numberwithin{figure}{section}
\numberwithin{table}{section}

Consider the scalar quasilinear parabolic differential equation
\begin{equation}\label{PDE}
    u_t = a(\theta,\phi,u,\nabla u)\Delta_{\mathbb{S}^2} u+f(\theta,\phi,u,\nabla u)
\end{equation}
with initial data $u(0,\theta,\phi)=u_0(\theta,\phi)$ such that $a,f\in C^2$, satisfy the strict parabolicity condition $a(\theta,\phi,u,\nabla u)\geq\epsilon>0$, and $\Delta_{\mathbb{S}^2}$ is the Laplace-Beltrami operator on the sphere $\mathbb{S}^2$. In coordinates, the angle variables are $(\theta,\phi)\in [0,\pi]\times [0,2\pi]$ with Neumann boundary condition in $\theta$ and periodic boundary in $\phi$.

Suppose that solutions $u(t,\theta,\phi)$ are \emph{axisymmetric}, that is, they are independent of rotations with respect to the angle $\phi$ and depend only in $\theta$. Hence, $u(t,\theta)$ solves the following equation
\begin{equation}\label{axisymPDE}
    u_t =  a(\theta,u,u_\theta)\left[u_{\theta\theta}+\frac{u_{\theta}}{\tan(\theta)}\right] +f(\theta,u,u_\theta)
\end{equation}
with initial data $u(0,\theta)=u_0(\theta)$, where $\theta \in [0,\pi]$ has Neumann boundary. Even though the equation has a singular coefficient at the boundaries $\theta=0$ or $\pi$, solutions are still regular. 

The equation \eqref{axisymPDE} defines a semiflow denoted by $(t,u_0)\mapsto u(t)$ in a Banach space $X^\alpha:=C^{2\alpha+\beta}([0,\pi])$. We suppose that $2\alpha+\beta>1$ so that solutions are at least $C^1$. The appropriate functional setting is described in Section \ref{sec:func}. 

In order to study the long time behavior of \eqref{axisymPDE}, we suppose that $f$ satisfies the following conditions
\begin{align}\label{diss}
    f(x,u,0) \cdot u&<0\nonumber\\
    |f(x,u,p)|&<f_1(u)+f_2(u)|p|^\gamma\nonumber\\
    \frac{|a_x|}{1+|p|}+|a_u|+|a_p|\cdot [1+|p|]&\leq f_3(|u|)\\
    0<\epsilon \leq a(x,u,p) &\leq \delta \nonumber 
\end{align}
where the first condition holds for $|u|$ large enough, uniformly in $x$, the second for all $(x,u,p)$ for continuous $f_1,f_2$ and $\gamma<2$, the third for continuous $f_3$ and $\epsilon,\delta>0$. 

Those conditions imply that $|u|$ and $|u_x|$ are bounded. Hence, the semiflow is dissipative: trajectories $u(t)$ eventually enter a large ball in the phase-space $X^\alpha$. See Chapter 6, Section 5 in \cite{Ladyzhenskaya68}. Also \cite{KruzhkovOleinek61} and \cite{BabinVishik92}. 

Moreover, these hypotheses guarantee that there exists a nonempty global attractor $\mathcal{A}$ of \eqref{axisymPDE}, which is the maximal compact invariant set. Equivalently, it is the set of bounded trajectories $u(t)$ in the phase-space $X^\alpha$ that exist for all $t\in \mathbb{R}$. See \cite{BabinVishik92}.

The goal of this paper is to decompose $\mathcal{A}$ into smaller invariant sets, and describe how those sets are related.

For the statement of the main theorem that describes the global attractor $\mathcal{A}$, denote by the \emph{zero number} $z(u_*)$ the number of strict sign changes of a continuous function $u_*(\theta)$. Recall that the \emph{Morse index} $i(u_*)$ of an equilibrium $u_*$ is given by the number of positive eigenvalues of the linearized operator at such equilibrium, that is, the dimension of the unstable manifold of $u_*$.

We say that two different equilibria $u_-,u_+$ of \eqref{axisymPDE} are \emph{adjacent} if there does not exist an equilibrium $u_*$ between $u_-$ and $u_+$ at $\theta=0$ satisfying 
\begin{equation}    
z(u_--u_*)=z(u_{+}-u_*).
\end{equation} 
 
This notion was firstly described by Wolfrum \cite{Wolfrum02}.

Both the zero number and Morse index can be computed from a permutation of the equilibria, as it was done in \cite{FuscoRocha} and \cite{FiedlerRocha96}. Such permutation is called the \emph{Sturm Permutation}. We construct an analogous permutation for the case of boundary singularity in Section \eqref{sec:perm}, as in \cite{FiedlerRocha96}. For such, it is required that the flow of the equilibria equation of \eqref{axisymPDE} exists for all $\theta\in [0,\pi]$. 


\begin{thm}\emph{\textbf{Sturm Attractor}} \label{attractorthm}

Consider $a,f\in C^2$ satisfying the growth conditions \eqref{diss}. Suppose that all equilibria for the equation \eqref{axisymPDE} are hyperbolic. Then, 
\begin{enumerate}
    \item the global attractor $\mathcal{A}$ of \eqref{axisymPDE} consists of finitely many equilibria $\mathcal{E}$ and their heteroclinic orbits $\mathcal{H}$.
    \item there exists a heteroclinic $u(t)\in\mathcal{H}$ between $u_-,u_+\in\mathcal{E}$ such that 
    \begin{equation*}
        u(t)\xrightarrow{t\to \pm\infty} u_{\pm}
    \end{equation*}
    if, and only if, $u_-$ and $u_+$ are adjacent and $i(u_-)>i(u_+)$.
\end{enumerate}
\end{thm}

The first claim follows due to the existence of a Lyapunov functional constructed by Matano \cite{Matano88} and Zelenyak \cite{Zelenyak68}. A modification of such functional for the case of singular coefficients is done in Section \ref{sec:func}.

The second claim answers the question of which equilibria connect to which other. This geometric description was carried out by Hale and do Nascimento \cite{HaleNascimento83} for the Chafee Infante problem, by Brunovsk\'y and Fiedler \cite{FiedlerBrunovsky89} for $f(u)$, by Fiedler and Rocha \cite{FiedlerRocha96} for $f(x,u,u_x)$, and for quasilinear equations by the author in \cite{Lappicyquasi17}. Such attractors are known as \emph{Sturm attractors}.

Constructing the Sturm attractor for the equation \eqref{axisymPDE} is problematic due to its singular coefficient. It is the aim of this paper to modify the existing theory for such boundary singularity and still obtain a Sturm attractor. 

In particular, we compute the attractor explicitely for the example of Chafee-Infante type nonlinearity with singular boundary coefficients. This attractor could be used as an application of the Einstein Hamiltonian equation, as in \cite{LappicyBlackHoles} and \cite{LappicyPhD}.

\begin{cor}\emph{\textbf{Chafee-Infante Attractor}} \label{singCIINTRO}

Consider $f(\lambda, \theta,u,u_\theta)=\lambda a(\theta,u,u_\theta)u(1-u^2)$ in the equation \eqref{axisymPDE}. Let $\lambda\in(\lambda_k,\lambda_{k+1})$, where $\lambda_k$ is the $k$-th eigenvalue of the axisymmetric Laplacian with $k\in\mathbb{N}_0$. 

Then, there are $2k+3$ hyperbolic equilibria $u_1,...u_{2k+3}$ and its attractor $\mathcal{A}$ is below in Figure \ref{FIGCOR}, where arrows denote heteroclinics.
\begin{figure}[ht]\centering
\begin{tikzpicture}
\filldraw [black] (0,0) circle (3pt) node[anchor=south]{$u_{k+2}\equiv 0$};
\filldraw [black] (-1,-1) circle (3pt) node[anchor=east]{$u_{k+1}$};
\filldraw [black] (1,-1) circle (3pt) node[anchor=west]{$u_{k+3}$};

\filldraw [black] (-1,-2.4) circle (0.5pt);
\filldraw [black] (-1,-2.5) circle (0.5pt);
\filldraw [black] (-1,-2.6) circle (0.5pt);
\filldraw [black] (1,-2.4) circle (0.5pt);
\filldraw [black] (1,-2.5) circle (0.5pt);
\filldraw [black] (1,-2.6) circle (0.5pt);

\filldraw [black] (-1,-4) circle (3pt) node[anchor=east]{$u_{2}$};
\filldraw [black] (1,-4) circle (3pt) node[anchor=west]{$u_{2k+2}$};
\filldraw [black] (-1,-5) circle (3pt) node[anchor=east]{$u_{1}\equiv -1$};
\filldraw [black] (1,-5) circle (3pt) node[anchor=west]{$u_{2k+3}\equiv +1$};

\draw[thick,->] (0,0) -- (-0.9,-0.9);
\draw[thick,->] (0,0) -- (0.9,-0.9);

\draw[thick,->] (-1,-1) -- (0.9,-1.9);
\draw[thick,->] (1,-1) -- (-0.9,-1.9);
\draw[thick,->] (-1,-1) -- (-1,-1.87);
\draw[thick,->] (1,-1) -- (1,-1.87);

\draw[thick,->] (-1,-3) -- (0.9,-3.9);
\draw[thick,->] (1,-3) -- (-0.9,-3.9);
\draw[thick,->] (-1,-3) -- (-1,-3.87);
\draw[thick,->] (1,-3) -- (1,-3.87);

\draw[thick,->] (-1,-4) -- (0.9,-4.9);
\draw[thick,->] (1,-4) -- (-0.9,-4.9);
\draw[thick,->] (-1,-4) -- (-1,-4.87);
\draw[thick,->] (1,-4) -- (1,-4.87);
\end{tikzpicture}
\caption{Global attractor $\mathcal{A}$ of Chafee-Infante type} \label{FIGCOR}
\end{figure}
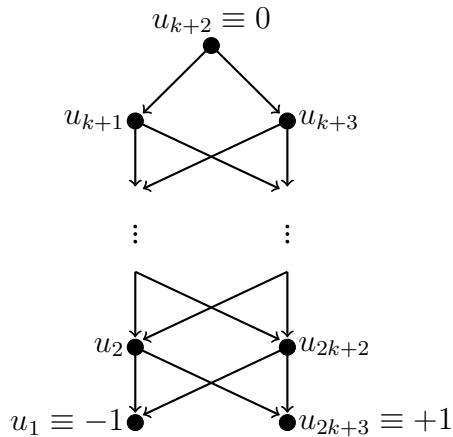
\end{cor}

This corollary is proved by constructing a Sturm permutation of the axisymmetric Chafee-Infante, yielding the same permutation as the usual Chafee-Infante problem. Hence, their attractors are geometrically (connection-wise) the same and their only difference lies in the form of the equilibria and the domain of the parameter $\lambda$.

The remaining sections are organized as follows. 

We firstly introduce the functional setting in Section \ref{sec:func}, and construct a Lyapunov functional for the singular case by modifying Matano's arguments from \cite{Matano78}, and its generalization for fully nonlinear equations \cite{LappicyFiedler18}. In particular this implies that the attractor consists of equilibria and heteroclinics. 

Then, we focus on the connection problem. All the necessary information about the adjacency, namely the zero numbers and Morse indices, are encoded in a permutation of the equilibria, which is described in Section \ref{sec:perm}. This was done firstly by \cite{FuscoRocha}, and here is modified for the singular case. 

In Section \ref{sec:globalsturm}, it is proven the dropping lemma for the singular case, as well as some consequences. This is a fundamental result for the attractor construction that dates back to Sturm and is done by modifying arguments of Chen and Pol\'acik \cite{ChenPolacik96}, where they proved such result for a singular coefficient at only one boundary value. Then all the previous tools are put together to construct the attractor in Section \ref{sec:globalsturm}, as it was done \cite{FiedlerRocha96}. 

Lastly, Section \ref{sec:CI} gives an example of the developed theory and constructs the attractor for the axisymmetric Chafee-Infante problem.

\section{Proof of main result}

\subsection{Functional setting} \label{sec:func}

The Banach space used on the upcoming theory consists on subspaces of H\"older continuous functions $C^\beta(\mathbb{S}^2)$ with $\beta\in (0,1)$. A more precise description is given below, following Lunardi \cite{Lunardi95} for the euclidean case, and Huang  \cite{Huang15} for the case on manifolds. The notation $C^{\beta}$ for some $\beta\in\mathbb{R}$ indicates that $\beta=[\beta]+\{\beta\}$, where the integer part $[\beta]\in\mathbb{N}$ denotes the $[\beta]$-times differentiable functions whose $[\beta]$-derivatives are $\{\beta\}$-Hölder, where $\{\beta\}\in[0,1)$ is the fractional part of $\beta$.

The equation \eqref{PDE} is seen as an abstract differential equation on a Banach space,
\begin{equation}\label{linequiv}
    u_t=Au+g(u)    
\end{equation}
where $A:D(A)\rightarrow \tilde{X}$ is the linearization of the right-hand side of \eqref{PDE} at the initial data $u_0(x)$, and the Nemitskii operator $g$ of the remaining nonlinear part, which takes values in $\tilde{X}$. The spaces considered are $\tilde{X}:=C^\beta(\mathbb{S}^2)$, and $D(A)=C^{2+\beta}(\mathbb{S}^2)\subset \tilde{X}$ is the domain of the operator $A$, where $\beta\in (0,1)$.

We consider the interpolation spaces $C^{2\alpha+\beta}(\mathbb{S}^2)$ between $D(A)$ and $\tilde{X}$ with $\alpha\in (0,1)$ such that $A$ generates a strongly continuous semigroup in $C^{2\alpha+\beta}(\mathbb{S}^2)$, and hence the equation \eqref{PDE} with the dissipative conditions \eqref{diss} defines a dissipative dynamical system in $\tilde{X}^\alpha$. We suppose that $2\alpha+\beta>1$ so that solutions are at least in $C^1(\mathbb{S}^2)$. Moreover, due to the Sobolev embedding, we know that $C^{2\alpha+\beta}(\mathbb{S}^2)\subseteq L^2(\mathbb{S}^2)$, and hence $C^{2\alpha+\beta}(\mathbb{S}^2)$ inherits an inner product, once its functions are considered as $L^2(\mathbb{S}^2)$ functions. Note all these spaces have metrics depending on the metric of the sphere.

In particular, it settles the theory of existence and uniqueness. For certain qualitative properties of solutions, such as the existence of invariant manifolds tangent to the linear eigenspaces, one needs to know the spectrum of $A$.

Now, we consider the restriction of the flow $u(t)$ in $C^{2\alpha+\beta}(\mathbb{S}^2)$ to the invariant subspace $X^\alpha:=C^{2\alpha+\beta}([0,\pi])$ which consists of functions that do not depend on the angle $\phi$. We now prove $X^\alpha$ is indeed invariant. Consider the projection $P:\tilde{X}^\alpha\to X^\alpha$. Let $w(t):=u(t)-Pu(t)$, where $Pu(t)$ is the restricted flow with the same initial data $Pu_0=u_0\in C^{2\alpha+\beta}([0,\pi])$. Note $w$ has initial data $w_0\equiv 0$ and satisfies a linear PDE. The maximum principle implies that $w(t)\equiv 0$, and hence the full flow $u(t)$ is equal to the restricted flow $Pu(t)$, which do not depend on $\phi$.

Hence, the equation \eqref{axisymPDE} can be rewritten as an equation in ${X^\alpha}$ like \eqref{linequiv} with the restricted operator $A|_{X^\alpha}$,
where the metric in $X^\alpha$ is the sum of the $C^{[2\alpha+\beta]}$-norm and the Hölder $C^{\{ 2\alpha+\beta \}}$-seminorm,
\begin{equation*} 
    ||u(\theta)||_{X^\alpha}=\sum_{k=0}^{[2\alpha+\beta]} \max_{\theta\in [0,\pi]} |\partial_\theta^k u(\theta)|+\max_{\theta_1,\theta_2\in [0,\pi]} \frac{|\partial^{[2\alpha+\beta]}_\theta u(\theta_1)-\partial^{[2\alpha+\beta]}_\theta u(\theta_2)|}{|\theta_1-\theta_2|_{arc}^{\{2\alpha+\beta\}}}    
\end{equation*}
where $|.|$ denotes the usual norm in $\mathbb{R}$, $\partial^{k}_\theta$ denotes the $k^{\text{th}}$-derivative with respect to $\theta$, and the distance in the axial-arc within the sphere is given by $|\theta_1-\theta_2|_{arc}=\int_{\theta_1}^{\theta_2}\sin(\theta) d\theta$. 

Moreover, $X^\alpha \subseteq L^2_w([0,\pi])$, and we also have an induced inner product in $X^\alpha$ given by
\begin{equation*}
        ||u(\theta)||_{L^2_w([0,\pi])}=\left(\int_{0}^\pi u^2(\theta) \sin(\theta) d\theta \right)^{\frac{1}{2}}
\end{equation*}
where the space $L^2_w([0,\pi])$ has weight $w:=\sin(\theta)$ that tames the singular term and is given by the usual spherical metric restricted on the axially symmetric arc within the sphere, parametrized by the angle $\theta\in [0,\pi]$. 

The operator $A$ is a self-adjoint singular Sturm-Liouville operator on the space $L^2_w([0,\pi])$. Its spectrum consists of real and simple eigenvalues $\lambda_k=k(k+1)$ for $k\in\mathbb{N}_0$ with Legendre polynomials $\phi_k= P_k(\cos{\theta})$ as corresponding eigenfunctions, which form an orthonormal basis of $L^2_w([0,\pi])$. Note all $\phi_k \in X^\alpha$, and hence are also a basis of $X^\alpha$. 


These yield the existence properties of invariant manifolds.

\begin{thm} \label{hierarchy}
    \emph{\textbf{Filtration of Invariant Manifolds} \cite{Mielke91}}
    
    Let $u_*$ be a hyperbolic equilibrium of \eqref{axisymPDE} with Morse index $n:=i(u_*)$. Then there exists a filtration of the \emph{unstable manifold}
    \begin{equation*}
        W^u_0(u_*)\subset ... \subset W^u_{n-1}(u_*)=W^u(u_*)
    \end{equation*}
    where each $W^u_k$ has dimension $k+1$ and tangent space at $u_*$ spanned by $\phi_0,...,\phi_{k}$. 
    
    Analogously, there is a filtration of the \emph{stable manifold}
    \begin{equation*}
        ... \subset W^s_{n+1}(u_*)\subset W^s_{n}(u_*)=W^s(u_*)
    \end{equation*}
    where each $W^s_k$ has codimension $k$ and tangent space at $u_*$ space spanned by $\phi_{k},\phi_{k+1},...$.
\end{thm}

Note that the above index labels are not in agreement with the dimension of each submanifold within the filtration, but it is with the number of zeros its corresponding eigenfunction has. For example, an eigenfunction $\phi_k$ corresponding to the eigenvalue $\lambda_k>0$ has $k$ simple zeroes, whereas the $\dim(W^u_k)=k+1$. 

An important property is the behavior of solutions within each submanifold of the above filtration of the unstable or stable manifolds.

\begin{thm} \label{ConvEF}
    \emph{\textbf{Linear Asymptotic behavior} \cite{Henry81},  \cite{Angenent86}, \cite{FiedlerBrunovsky86}} 
    
    Consider a hyperbolic equilibrium $u_*$ with Morse index $n:=i(u_*)$ and a trajectory $u(t)$ of \eqref{axisymPDE}. The following holds,
    \begin{enumerate}
    \item If $u(t)\in W^u_k(u_*)\backslash W^u_{k-1}(u_*)$ with $k=0,...,i(u_*)-1$, then
    \begin{equation*}
        \frac{u(t)-u_*}{||u(t)-u_*||}\xrightarrow{t\rightarrow -\infty} \pm \phi_k.
        \end{equation*}
    \item If $u(t)$ in $W^s_k(u_*)\backslash W^s_{k+1}(u_*)$ with $k\geq i(u_*)$, then
        \begin{equation*}
        \frac{u(t)-u_*}{||u(t)-u_*||}\xrightarrow{t\rightarrow \infty} \pm \phi_k..
        \end{equation*}
    \end{enumerate}

    where the convergence takes place in $X^\alpha\subseteq C^1$, and $W^u_{-1}(u_*)=\emptyset$. 
    
    The conclusions of 1. and 2. also hold true by replacing the difference $u(t)-u_*$ with the tangent vector $u_t$.
\end{thm}

The reason this theorem works for both the tangent vector $v:=u_t$ or the difference $v:=u_1-u_2$ of any two solutions $u_1$ and $u_2$ of the nonlinear equation \eqref{axisymPDE} is because they satisfy a linear equation of the type
\begin{equation}\label{linaxisymPDE}
    v_t=a(t,\theta)\left[v_{\theta\theta}+\frac{v_\theta}{\tan(\theta)}\right]+b(t,\theta)v_\theta+c(t,\theta)v
\end{equation}
where $\theta\in(0,\pi)$ has Neumann boundary conditions, $a(t,\theta),b(t,\theta)$ and $c(t,\theta)$ are bounded.



Now we show that there exists a Lyapunov function, as it was done by Zelenyak \cite{Zelenyak68} and Matano \cite{Matano82}. We modify Matano's construction bearing in mind that the metric on the sphere induces a space with weighted norms, and this weight should be incorporated into the construction of the Lyapunov function. Therefore, the flow in $X^\alpha$ is gradient with respect to the induced inner product of $L^2_w$. As a consequence of the Lyapunov function, bounded trajectories tend to equilibria.

\begin{lem}
    \emph{\textbf{Lyapunov Function}}
    
    There exists a Lagrange function $L$ such that 
    \begin{equation}\label{Lyap}
        E:= \int_{0}^\pi L(\theta,u,u_\theta) \sin(\theta) d\theta
    \end{equation}
    is a Lyapunov function for the equation \eqref{axisymPDE}.
\end{lem}

Note that in the case that the nonlinearity $f$ does not depend on $u_\theta$, then the Lagrange functional $L(\theta,u,u_\theta):=\frac{1}{2}u_\theta^2-F(\theta,u)$ yields a Lyapunov function $E$, where $F$ is the primitive function of $f$. 
Indeed,
\begin{equation*}
    \frac{dE}{dt}=-\int_{0}^\pi |u_t|^2 \sin(\theta) d\theta \leq 0. 
\end{equation*}

For nonlinearities of the type $f(\theta,u,u_\theta)$, we obtain a Lyapunov function such that
\begin{equation}\label{imlisteningtocoltrane}
    \frac{dE}{dt}:= -\int_{0}^\pi  |u_t|^2 \frac{L_{pp}\sin(\theta)}{a}   d\theta\leq 0
\end{equation}
where $p:=u_\theta$ and $L$ satisfy the convexity condition $L_{pp}>0$. Hence, the case that $f$ does not depend on $u_\theta$ is seen as a particular case when $L_{pp}=1$.
\begin{pf}
Let $p:=u_\theta$ and differentiate \eqref{Lyap} with respect to $t$,
\begin{equation*}
    \frac{dE}{dt}= \int_{0}^\pi \left[ L_u u_t + L_p u_{\theta t}\right] \sin(\theta) d\theta .
\end{equation*}    

Integrating the second term by parts and noticing that the $\sin(\theta)$ is 0 at the boundaries,
\begin{align*}
    \frac{dE}{dt}&= \int_{0}^\pi \left[ L_u \sin(\theta)-\frac{d}{d\theta}(L_p \sin(\theta)) \right] u_t d\theta\\
    &= \int_{0}^\pi \left[ (L_u -L_{p\theta} -L_{pu}u_\theta -L_{pp}u_{\theta\theta})\sin(\theta) -L_{p}\cos{\theta}\right] u_t d\theta.
\end{align*}

Substitute \eqref{axisymPDE} casted as $u_{\theta\theta} \sin(\theta)=[u_t \sin(\theta)-f\sin(\theta)]/a-u_\theta \cos(\theta)$,
\begin{align*}
    \frac{dE}{dt}&=\int_{0}^\pi (L_u-L_{p\theta}-L_{pu}u_\theta +L_{pp}\frac{f}{a})\sin(\theta)u_td\theta \\
    & \ \ \ \  +\int_{0}^\pi(L_{pp}u_\theta-L_{p})\cos(\theta) u_t d\theta -\int_{0}^\pi \frac{L_{pp}}{a} u_t^2 \sin(\theta) d\theta.
\end{align*}

To obtain \eqref{imlisteningtocoltrane}, we guarantee that there exists a function $L$ satisfying
\begin{equation}\label{Lyap0}
\left(L_u-L_{p\theta}-L_{pu}p +L_{pp}\frac{f}{a}\right)\sin(\theta) +(L_{pp}p-L_{p})\cos(\theta)=0
\end{equation}
for all $u,p\in\mathbb{R}$ and $\theta\in[0,\pi]$. 

Differentiating this equation with respect to $p$, some of the terms cancel, yielding
\begin{equation}\label{Lyap0_p}
\left[-L_{pp\theta}-L_{ppu}p+L_{ppp}\frac{f}{a}+L_{pp}\left(\frac{f}{a}\right)_p\right]\sin(\theta) +(L_{ppp}p)\cos(\theta)=0.
\end{equation}

To make sure that $L_{pp}>0$, introduce $g=g(\theta,u,p)$ through $L_{pp}=\exp(g)>0$. Hence, $g$ satisfies the following linear first order differential equation,
\begin{equation}\label{Lyapg}
    \left[\left(g_{\theta}+g_{u}p-g_{p}\frac{f}{a}-\left(\frac{f}{a}\right)_p\right)\sin(\theta)-(g_{p}p)\cos(\theta)\right] \exp(g)=0.
\end{equation}

Or equivalenty,
\begin{equation*}
    \left[g_{\theta}+g_{u}p+g_{p}\left(-\frac{f}{a}-p\frac{\cos(\theta)}{\sin(\theta)} \right) \right]\sin(\theta)=\left(\frac{f}{a}\right)_p\sin(\theta).    
\end{equation*}

This can be solved through the method of characteristics: along the solutions of 
\begin{align*}
\begin{cases}
    {u_\theta}&=p\\
    {p_\theta}&=-\frac{f}{a}-p \frac{\cos(\theta)}{\sin(\theta)}
\end{cases}
\end{align*}
the function $g$ must satisfy
\begin{equation*}
\frac{dg}{d\theta}=\left(\frac{f}{a}\right)_p.
\end{equation*}


Note that characteristics solve the equation for equilibria. If solutions of such equations exist for all initial conditions $(u,p)\in \mathbb{R}^2$ at $\theta=0$, and all $\theta\in[0,\pi]$, 
we obtain a global solution $g$ of \eqref{Lyapg} with some initial data, for example, $g(0,u,p)\equiv 0$. 

It is still needed to ascend from a function $g$ satisfying \eqref{Lyapg} to a function $L$ satisfying \eqref{Lyap0}. A choice for $L$ such that $L_{pp}=\exp(g)$ can be obtained by integrating this relation twice with respect to $p$, yielding a solution of \eqref{Lyap0_p},
\begin{equation*}
    L(\theta,u,p):=\int_0^p \int_0^{p_1} \exp(g(\theta,u,p_2))dp_2dp_1 + G(\theta,u).
\end{equation*}

To show that such $L$ is also a a solution of \eqref{Lyap0}, we have to restrict which $G$ are allowed. Recall that \eqref{Lyap0_p} was obtained through differentiating \eqref{Lyap0} with respect to $p$. That means that the left-hand side of \eqref{Lyap0} is independent of $p$, since it is equal to $0$. Hence it is satisfied for all $p$, if it holds for $p=0$. 

At $p=0$, the construction of $L$ yields that $L_p=L_{p\theta}=0$ and $L_u=G_u$. Plugging it in the equation \eqref{Lyap0} at $p=0$, it yields $(G_u+L_{pp}f/a)\sin(\theta)=0$. Hence, $G_u+L_{pp}f/a=0$, that is, $G_u=-\exp(g)f/a$. Integrating in $u$,
\begin{equation*}
    G(\theta,u):=- \int_0^{u} \frac{f(\theta,u_1,0)}{a(\theta,u_1,0)} exp(g(\theta,u_1,0))du_1 
\end{equation*}
\begin{flushright}
	$\blacksquare$
\end{flushright}
\end{pf}

Note that one can do a similar construction of a Lyapunov function without assuming that the $\sin(\theta)$ appears in the integrand, as in \eqref{Lyap}. But such coefficient will appear once the differential equation is plugged in the Ansatz for the Lyapunov functional.

Moreover, Matano's construction can be adapted to more general singular Sturm Liouville operators of the form $\frac{\partial_\theta(r(\theta) \partial_\theta)}{w(\theta)}$, if the weight $\sin(\theta)$ within the integrand is replaced by $w(\theta)$. Hence, the Lyapunov function will decay in the $L^2_w$ norm with appropriate weighted metric $w$.

Therefore, the LaSalle invariance principle holds and implies that bounded solutions converge to equilibria, and any $\alpha,\omega$-limit set consist of a single equilibrium. See \cite{Matano88}. Moreover, due to hyperbolicity, equilibria are isolated and due to dissipativity, there are finitely many of them. Hence, the global attractor consists of finitely many equilibria, and their heteroclinic connections, yielding the first part of the main result. See \cite{Henry81} and \cite{BabinVishik92}.

\subsection{Sturm permutation}\label{sec:perm}

The next step on our quest to find the Sturm attractor is to construct a permutation associated to the equilibria, which is done using shooting methods. This enables the computation of the Morse indices and zero number of equilibria. That was firstly done by Fusco and Rocha \cite{FuscoRocha} using methods also described by Fusco, Hale and Rocha in \cite{Rocha85}, \cite{RochaHale85}, \cite{rocha88}, \cite{Rocha94} and \cite{FuscoHale85}.

The equilibria equation associated to \eqref{axisymPDE} can be rewritten as
\begin{equation}\label{PDEeq}
    0=\frac{a(\theta,u,u_\theta)}{\sin(\theta)}\frac{d}{d\theta}[u_{\theta}\sin(\theta)]+f(\theta,u,u_\theta)
\end{equation}
for $\theta \in [0,\pi]$ with Neumann boundary conditions and $a>\epsilon>0$. 

In order to get rid of the singularities at $\theta=0$ and $\pi$, rescale the system by $\tau(\theta):=\ln(\tan(\theta/2))\in(-\infty,\infty)$, which maps the singularities at $\theta=0,\pi$ to $\tau=\pm \infty$. Also, add the equation $\theta_\tau=\sin(\theta)$ to obtain an autonomous system,
\begin{align*}
\begin{cases}
    0&=a(\theta,u,u_\tau /\sin(\theta))u_{\tau\tau}+ f(\theta,u,u_\tau /\sin(\theta))\sin^2(\theta)\\
	\theta_\tau&=\sin(\theta).    
\end{cases}
\end{align*}

Moreover, reduce the equation to first order system through $p:=u_\tau$. Hence, 
\begin{align}\label{shootflow}
\begin{cases}
    u_{\tau}&= p\\
	p_{\tau}&=- \frac{f\left(\theta,u,p/\sin(\theta)\right)}{a(\theta,u,p/\sin(\theta))}\sin^2(\theta)\\
	\theta_\tau&=\sin(\theta)
\end{cases}
\end{align}
where the Neumann boundary condition becomes 
$\lim_{\tau\rightarrow\pm\infty} p(\tau)=0$, since the Neumann boundary in changed of coordinates yields 
\begin{equation}\label{shootBC}
    0=\displaystyle \lim_{\theta\rightarrow 0,\pi} u_\theta=\lim_{\tau\rightarrow\pm\infty} u_\tau\frac{d\tau}{d\theta}=\lim_{\tau\rightarrow\pm\infty} p(\tau)\cosh(\tau)
\end{equation}
and $\lim_{\tau\rightarrow\pm\infty}\cosh(\tau)\to\infty$. This forces exponential decay of $p$.

Note that the term $\sin^2(\theta)$ cuts off the reaction $f$, being $1$ at the equator and decaying to $0$ near the poles. This means that the diffusion near the poles is stronger. 

In the nonsingular case, the idea to find equilibria \eqref{axisymPDE} is as follows. They must lie in the line 
\begin{equation*}
    L_0:=\{(\theta,u,p)\in\mathbb{R}^3 \textbf{ $|$ }  (\theta,u,p)=(0,d,0) \text{ and } d\in\mathbb{R}\}  
\end{equation*}
due to Neumann boundary at $\theta=0$. Then, evolve this line under the flow of the equilibria differential equation and intersect it with an analogous line $L_\pi$ at $\theta=\pi$, so that it also satisfies Neumann at $\theta=\pi$. This reasoning does not work for the singular case, since $L_0$ is a line of equilibria and is invariant under the shooting flow \eqref{shootflow}. A new approach is needed.

In the singular case, the linearization of \eqref{shootflow} at each point in $L_0$ has eigenvalues $\lambda_1=1$ and $\lambda_2=\lambda_3=0$ with respective generalized eigenvectors $v_1=(0,0,1),v_2=(1,0,0),v_3=(0,1,0)$. Hence, there is an one dimensional unstable direction given by the $\theta$-axis, and two center directions given by the invariant plane $\{(u,p,0)\in\mathbb{R}^3\}$. 

Furthermore, each point $(0,d,0)\in L_0$ has an one dimensional strong unstable manifold $W^u(0,d,0)$, which is locally a graph $\{(\theta,u^u(\theta,d),p^u(\theta,d))\in\mathbb{R}^3\} $. See \cite{GuckenheimerHolmes83}. 
The collection of all these strong unstable manifolds defines the \emph{unstable shooting manifold} $M^u$,
\begin{equation*}
    M^u:= \bigcup_{d\in\mathbb{R}} W^u(0,d,0).
\end{equation*}

Similarly, each point $(0,e,0)\in L_\pi$ has an one-dimensional strong stable manifold given locally by the graph $\{(\theta,u^s(\theta,e),p^s(\theta,e))\in\mathbb{R}^3\}$, and its collection defines the \emph{stable shooting manifold} $M^s$,
\begin{equation*}
    M^s:= \bigcup_{e\in\mathbb{R}} W^s(0,e,0).
\end{equation*}


We assume that solutions of \eqref{shootflow} are defined for all $\theta\in[0,\pi]$ and any initial data $(u,p)$. Hence,
the shooting manifolds will exist globally and for any initial data. 

Denote by $M^u_\theta$ the cross-section of $M^u$ for some fixed $\theta\in [0,\pi]$. This is a curve parametrized by $d\in\mathbb{R}$. Similarly, $M^s_\theta$ is a curve parametrized by $e\in\mathbb{R}$. 

We obtain the following characterization of equilibria its Morse indices and zero numbers, through the shooting manifolds, similar to \cite{Rocha85} and \cite{Hale99}.

\begin{lem}
    \textbf{Equilibria Through Shooting}
    \begin{enumerate}
    \item The set of equilibria $\mathcal{E}$ of \eqref{axisymPDE} is in one-to-one correspondence with $M^u_{\theta}\cap M^s_{\theta}$ for any $\theta\in[0,\pi]$.
    \item An equilibrium point corresponding to fixed $d\in\mathbb{R}$ and $e\in\mathbb{R}$ is hyperbolic if, and only if, $W^u(0,d,0)$ intersects $W^s(0,e,0)$ transversely.
    \item  If $u_*$ correspond to a hyperbolic equilibrium of \eqref{axisymPDE}, then its Morse index is given by $i(u_*)=1+\lfloor\frac{\zeta(\theta_0)}{\pi}\rfloor$ where $\zeta(\theta_0)$ is the angle between $M^u$ and $M^s$ measured clockwise at their intersection point $\theta_0$, and $\lfloor.\rfloor$ denotes the floor function. 
    \end{enumerate}
\end{lem}
\begin{pf}
To prove 1), note that a point in $M^u_{\theta}\cap M^s_{\theta}$ satisfies the equilibria equation by definition of the shooting manifolds. Moreover, the Neumann boundary conditions are also satisfied since solutions are in the appropriate stable/unstable manifolds.

Conversely, consider an equilibrium of \eqref{axisymPDE}. It must satisfiy the Neumann boundary conditions \eqref{shootBC}, which requires exponential convergence rate to $0$. This implies that the equilibrium must be both in the strong unstable $M^u$ and strong stable $M^s$ manifolds. Moreover, such manifolds intersect for some $\theta\in[0,\pi]$, because the equilibrium is continuous. By uniqueness and invariance of the shooting manifolds, they must also intersect for all $\theta\in[0,\pi]$.

Due to the uniqueness of the shooting differential equation \eqref{shootflow}, such correspondence above is one-to-one. 

To prove 2), consider an equilibrium $u_*$ corresponding to $d,e\in\mathbb{R}$. We compare the eigenvalue problem for $u_*$ and the differential equation satisfied by the angle of the tangent vectors of the shooting manifold. 

Introducing the $\tau$ variable, the eigenvalue problem for $u_*$ is obtained by linearizing the right hand side of the equation in order to obtain a linear operator, yielding
\begin{align*}
\begin{cases}
    \lambda u \sin^2(\theta)&=a_*u_{\tau\tau}+b_*u+c_*u_\tau \\
    \theta_\tau&=\sin(\theta)
\end{cases}
\end{align*}
with boundary conditions $\lim_{\tau\rightarrow\pm\infty} u_\tau(\tau)=0$, where 
\begin{align*}
    a_*(\theta)&:=a(\theta,u_*,\partial_\tau u_*)\\ b_*(\theta)&:=a_u(\theta,u_*,\partial_\tau u_*).(u_*)_{\theta\theta}+D_uf(\theta,u_*,\partial_\tau u_*)\sin^2(\theta)\\ c_*(\theta)&:=a_p(\theta,u_*,\partial_\tau u_*).(u_*)_{\theta\theta}+D_pf(\theta,u_*,\partial_\tau u_*)\sin^2(\theta).
\end{align*}

Rewriting the above system as a system of first order by $p:=u_\tau$,
\begin{align*}
\begin{cases}
    u_\tau&=p\\
    p_\tau&=-\frac{b_*u+c_*p-\lambda u\sin^2(\theta)}{a_*} \\
    \theta_\tau&=\sin(\theta)
\end{cases}
\end{align*}
with boundary conditions $ \lim_{\tau\rightarrow\pm\infty} p(\tau)=0$.

In polar coordinates $(u,p)=:(r\cos(\mu),-r\sin(\mu))$, the angle $\mu:=\arctan(\frac{p}{u})$ satisfies 
\begin{align}\label{EVpolar}
\begin{cases}
    \mu_\tau&= \sin^2(\mu)+ \frac{b_*u+c_*p-\lambda u\sin^2(\theta)}{a_*} \cos^2(\mu)\\
	\theta_\tau&=\sin(\theta)
\end{cases}
\end{align}
with $ \lim_{\tau=-\infty} \mu(\tau)=0$ and $ \lim_{\tau\rightarrow \infty} \mu(\tau)=k\pi$ for some $k\geq0$. 


On the other hand, $M^u_\theta$ is parametrized by $d\in\mathbb{R}$ and its tangent vector $(\frac{\partial u(\theta,d)}{\partial d},\frac{\partial p(\theta,d)}{\partial d})$ satisfies the following linearized equation, 
\begin{align}\label{tangentshoot}
\begin{cases}
    (u_{d})_{\tau}&= p_{d}\\
	(p_{d})_{\tau}&=-\frac{b^uu_d+c^up_d-\lambda u\sin^2(\theta)}{a^u}\\
	\theta_\tau&=\sin(\theta)
\end{cases}
\end{align}
with initial data $ \lim_{\tau\rightarrow -\infty} (u_{d},p_{d})=(1,0)$. Note that the linearization is considered along the unstable manifold given by the graph $\{(\theta,u^u(\theta),p^u(\theta))\in\mathbb{R}^3\}$, and the definition of $a^u,b^u,c^u$ are the same as $a_*,b_*,c_*$, except they are evaluated in the unstable manifold, instead of the equilibrium $u_*$.

In polar coordinates $(u_{d},p_{d})=:({\rho}\cos({\nu}),-{\rho}\sin({\nu}))$, where ${\nu}$ is the clockwise angle of the tangent vector of $M^u_\theta$ with the $u$-axis, 
\begin{align}\label{tangentpolar}
\begin{cases}
    {\nu}_\tau&= \sin^2({\nu})+ \frac{b^uu_d+c^up_d-\lambda u\sin^2(\theta)}{a^u} \cos^2({\nu})\\
	\theta_\tau&=\sin(\theta)
\end{cases}
\end{align}
with initial data $ \lim_{\tau\rightarrow -\infty} {\nu}(\tau,d)=0$. 

Similarly, the angle $\tilde{\nu}$ of the tangent vector of $M^s_\theta$ with the $u$-axis satisfies the equation \eqref{tangentpolar}, but with initial data $ \lim_{\tau\rightarrow \infty} \nu(\tau,e)=0$. 

Note that the equation \eqref{tangentpolar} that both angles $\nu$ and $\tilde{\nu}$ of the tangent vector satisfy is the same equation as the eigenvalue problem in polar coordinates \eqref{EVpolar} with $\lambda=0$, where each $\nu$ or $\tilde{\nu}$ encodes the boundary condition at $\tau=-\infty$ of $\infty$.

By hypothesis, the equilibrium $u_*$ corresponds to the pair of initial data $d,e\in\mathbb{R}$. That means that $M^u_{\theta_0}$ intersects $M^s_{\theta_0}$ for some fixed $\theta_0\in[0,\pi]$. 

Suppose that $u_*$ is not hyperbolic, that is, $\lim_{\tau\rightarrow \infty} \mu(\tau)=k\pi$ for $\lambda=0$ and some $k\in\mathbb{N}$. We compare this value with
the angle between the shooting curves at $\theta_0$. More precisely, it is proven that
\begin{equation}\label{totangle}
    \lim_{\tau\rightarrow \infty} \mu(\tau)=\nu(\theta_0)-\tilde{\nu}(\theta_0).
\end{equation} 

Indeed, for $\theta\in[0,\theta_0]$ the equations \eqref{EVpolar} and \eqref{tangentpolar} are the same, since both of them are linearized at the same orbit $u_*$, which corresponds to the unstable manifold of $(0,d,0)\in\mathbb{R}^3$. Since both of them have the same initial data, uniqueness implies
\begin{equation*}
    \mu(\theta_0)=\nu(\theta_0).
\end{equation*}

To obtain a relation between $\mu$ and $\tilde{\nu}$, consider the change of coordinates in the eigenvalue problem \eqref{EVpolar} as $\tilde{\mu}:=\mu-k\pi$. The equation $\eqref{EVpolar}$ is invariant under this transformation, since $\sin^2(\tilde{\mu}+k\pi)=\sin^2(\tilde{\mu})$. But the boundary condition changes at $\theta=\pi$, namely, $ \lim_{\tau\rightarrow \infty} \tilde{\mu}(\tau)=0$. Therefore, $\tilde{\mu}$ satisfies the same equation as the angle $\tilde{\nu}$, for $\theta\in[\theta_0,\pi]$. Hence, by uniqueness,
\begin{equation*}
    \mu(\theta_0)-k\pi=\tilde{\mu}(\theta_0)=\tilde{\nu}(\theta_0).
\end{equation*}

Subtracting these last two equations yields $k\pi=\nu(\theta_0)-\tilde{\nu}(\theta_0)$, that is, the intersection of the shooting manifolds is not transverse at their intersection point $\theta_0$.

Conversely, if the shooting manifolds are not transverse at some intersection point for $\theta_0$, then $k\pi=\nu(\theta_0)-\tilde{\nu}(\theta_0)$. 

Concatenate the solution $\nu$ from $M^u$ for $\theta\in[0,\theta_0]$ and initial data $\lim_{\tau=-\infty} \nu(\tau)=0$, together with $\tilde{\nu}$ from $M^s$ for $\theta\in [\theta_0,\pi]$ and initial data $\tilde{\nu}(\theta_0)=\nu(\theta_0)-k\pi$. Hence, the previous boundary conditions $\lim_{\tau=\infty} \tilde{\nu}(\tau)=0$ implies that $\lim_{\tau=\infty} \tilde{\nu}(\tau)=k\pi$, by considering the new initial data at $\theta=\theta_0$. Note such concatenated solution satisfy the equation \eqref{EVpolar} for the angle $\mu$ of the eigenvalue problem with $\lambda=0$. This implies there exists a solution $\mu$ of \eqref{EVpolar} and hence $\lambda=0$ is an eigenvalue. Thus, the equilibrium $u_*$ is not hyperbolic.

To prove 3), 
consider the solution $\mu(\tau,\lambda)$ of the eigenvalue problem in polar coordinates \eqref{EVpolar}. The Sturm oscillation theorem implies that 
\begin{equation}\label{anglelambda}
    \psi(\lambda):=\lim_{\tau\rightarrow \infty} \mu(\tau,\lambda)
\end{equation} 
is decreasing so that $\lim_{\lambda\rightarrow -\infty} \psi(\lambda)=\infty$ and $\lim_{\lambda\rightarrow \infty} \psi(\lambda)=-\pi/2$. Hence, there exists a decreasing sequence $\{ \lambda_k \}_{k\in N}$ to $-\infty$ such that $\psi(\lambda_k)=k\pi$ for $k\in \mathbb{N}$. This implies that there exists a solution of \eqref{EVpolar} for each $\lambda_k$ such that $\psi(\lambda_k)=k\pi$, and hence $\{ \lambda_k \}_{k\in N}$ are the eigenvalues.

Recall that the Morse index $i(u_*)$ is the number of positive eigenvalues of the linearization at $u_*$, that is 
\begin{equation*}
    ...<\lambda_{i(u_*)}<0<\lambda_{i(u_*)-1}<...<\lambda_0. 
\end{equation*}

Since $\psi(\lambda)$ is decreasing and $\lambda_{i(u_*)}$ are eigenvalues, then
\begin{equation*}
    i(u_*)\pi=\psi(\lambda_{i(u_*)})> \psi(0)> \psi(\lambda_{i(u_*)-1})=(i(u_*)-1)\pi .
\end{equation*}

Divide the above by $\pi$ and consider the integer value, yielding that $i(u_*)=\lfloor \frac{\psi(0)}{\pi}\rfloor +1$. It was noted in \eqref{totangle} that $\psi(0)=\nu(\theta_0)-\tilde{\nu}(\theta_0)$, which is exactly the angle between $M^u$ and $M^s$.
\begin{flushright}
	$\blacksquare$
\end{flushright}
\end{pf}

Hence, one can obtain a \emph{Sturm permutation} $\sigma$ by labeling the intersection points $u_i\in M^u_{\frac{\pi}{2}}\cap M^s_{\frac{\pi}{2}}$ firstly along $M^u_{\frac{\pi}{2}}$ following its para\-metrization given by $(\frac{\pi}{2},u^u(\frac{\pi}{2},d),p^u(\frac{\pi}{2},d))$ as $d$ goes from $-\infty$ to $\infty$. Namely,
\begin{equation*}
    u_1 < ... < u_N
\end{equation*} 
where $N$ denotes the number of equilibria. Secondly, label the intersection points along $M^s_{\frac{\pi}{2}}$ following its parametrization by $e\in\mathbb{R}$,
\begin{equation*}
    u_{\sigma(1)} < ... < u_{\sigma(N)}
\end{equation*} 

The Morse indices of equilibria and the zero number of difference of equilibria can be calculated through the Sturm permutation $\sigma$, as in \cite{Rocha91} and \cite{FiedlerRocha96}. This yields all necessary information for adjacency. The main tool for such proofs is the third part of the above Lemma: the rotation along the shooting curve increases the Morse index.




\subsection{Dropping lemma}\label{sec:drop}

Let the \emph{zero number} $z^t(u)$ count the number of strict sign changes in $\theta$ of a $C^1$ function $u(t,\theta)\not \equiv 0$, for each fixed $t$. More precisely,
  \[
    z^t(u):= \sup_k \left\{ 
                \begin{array}{ll}
                  \exists \text{ partition }\{ \theta_j\}_{j=1}^{k} \text{ of } [0,\pi]
                  \text{ such that }\\
                  u(t,\theta_j)u(t,\theta_{j+1})<0 \text{ for all } j=1,...,k
                \end{array}
          \right\}.
  \]          
and $z^t(u)=-1$ if $u\equiv 0$. In case $u$ does not depend on $t$, we simply write $z^t(u)=z(u)$.

A point $(t_0,\theta_0)\in\mathbb{R}\times [0,\pi]$ such that $u(t_0,\theta_0)=0$ is said to be a \emph{simple zero} if $u_\theta(t_0,\theta_0)\neq 0$ and a \emph{multiple zero} if $u_\theta(t_0,\theta_0)=0$.

The following result shows that the zero number of certain solutions of \eqref{axisymPDE} is nonincreasing in time $t$, and decreases whenever a multiple zero occur. Different versions of this well known fact are due to Sturm \cite{Sturm}, Matano \cite{Matano82}, Angenent \cite{Angenent88} and others. 
\begin{lem} \label{droplem}
    \emph{\textbf{Dropping Lemma}}
    
    Consider $v\not \equiv 0$ a solution of the linear equation \eqref{linaxisymPDE} for $t\in [0,T)$. Then, its zero number $z^t(v)$ satisfies
    \begin{enumerate}
    \item $z^t(v)<\infty$ for any $t\in (0,T)$.
    \item $z^t(v)$ is nonincreasing in time $t$.
    \item $z^t(v)$ decreases at multiple zeros $(t_0,\theta_0)$ of $v$, that is, 
    \begin{equation*}
        z^{t_0-\epsilon}(v)>z^{t_0+\epsilon}(v)
    \end{equation*} for any  sufficiently small $\epsilon>0$.
    \end{enumerate}
\end{lem}

Recall that both the tangent vector $u_t$ and the difference $u_1-u_2$ of two solutions $u_1,u_2$ of the nonlinear equation \eqref{axisymPDE} satisfy a linear equation as \eqref{linaxisymPDE}. Hence, the dropping lemma deals with the zero number of such solutions.

Below we give two different proofs. The first is an adaptation of Chen and Pol\'a\v{c}ik \cite{ChenPolacik96}, where the dropping lemma was proved for the case of a singular coefficient at one boundary point. The second by Angenent \cite{Angenent88}, where this lemma was proved for the case of regular coefficients. We also note that it is also possible to adapt the Newton polygon method done in Angenent \cite{Angenent90} and Angenent with Fiedler \cite{AngenentFiedler88}, but this is not pursued here, since this assumes that $a,f$ are analytic. 
\subsubsection{Proof 1}

This proof adapts Chen and Pol\'a\v{c}ik \cite{ChenPolacik96}. We cut off solutions nearby each boundary point so that it satisfies a differential equation with only one boundary singularity, and then apply the dropping lemma for such equations as it was proved in \cite{ChenPolacik96}.

We say two functions $u(t,\theta)$ and $v(t,\theta)$ have the \emph{same type of zeros} if for each fixed $t$, their zeros in $\theta$ coincide, together with their property of being simple or multiple. Mathematically, $u(t,\theta_0)=0$ if, and only if $v(t,\theta_0)=0$, for fixed $t$. Moreover, consider a zero $\theta_0$ of $u$ and $v$ for fixed $t$, then $u_\theta(t,\theta_0)=0$ if, and only if $v_\theta(t,\theta_0)=0$, .
\begin{lem}\label{adaptlem}
Suppose $u\not\equiv 0$ is a solution of \eqref{linaxisymPDE}. Then, there exists bounded functions v and d on $[t_1,t_2]\times[0,\pi]$ satisfying
\begin{equation}\label{cutoffPDE}
    v_t = v_{\theta\theta}+\frac{v_{\theta}}{\theta} + d(t,\theta)v
\end{equation}
where $\theta\in (0,\pi)$ has Neumann boundary conditions. Moreover, for a fixed  $\theta_1\in(0,\pi)$, the functions $u$ and $v$ have the same type of zeros for $\theta\in[0,\theta_1]$, whereas $v\neq 0$ for all $\theta\in[\theta_1,\pi]$.
\end{lem}

\begin{pf}
The idea is to localize the solution $u(t,\theta)$ for each $t$ and $\theta$ near the boundary $\theta=0$, and cut off whatever is far from it. Vaguely, this defines $v(t,\theta)$, and $d(t,\theta)$ is chosen accordingly so that one obtains the desired equation \eqref{cutoffPDE}.

Since the solution $u\not\equiv 0$, choose a point $\theta_1$ such that the solution is not zero at $\theta_1$ for a nonempty small interval of time $[t_1,t_2]$, by continuity in $t$. Moreover, due to continuity in $\theta$, choose $\theta_2\in(\theta_1,\pi)$ such that $u(t,\theta)\neq 0$ for $[t_1,t_2]\times [\theta_1,\theta_2]$. Without loss of generality, suppose that $u$ is positive for $[t_1,t_2]\times [\theta_1,\theta_2]$. Otherwise, consider $u(t,\theta)\mapsto -u(t,\theta)$.

Expand the singular term in power series as $\frac{1}{\tan(\theta)}=\frac{1}{\theta}+b(\theta)$, where $b(\theta)=\sum_{n=0}^\infty b_n \theta^{2n+1}$ is analytic in $\theta\in[0,\pi)$ and its coefficients $b_n$ are related to the Bernoulli numbers. Plugging this in \eqref{linaxisymPDE}, yields
\begin{equation*}
    u_t = u_{\theta\theta}+\frac{u_{\theta}}{\theta}+ b(\theta)u_{\theta}+ c(t,\theta)u.
\end{equation*}

Since $b(\theta)$ converges for $\theta\in [0,\pi)$ but not for $\theta=\pi$, this is how the singularity at $\theta=\pi$ is encoded in the new equation.

In order to get rid of $b(\theta)$, rescale the solution for $\theta\in[0,\theta_2]$ by $\tilde{u}(t,\theta):=\exp{(\frac{1}{2}\int_0^\theta b(y)dy)}u(t,\theta)$. Note $u(t,\theta)$ and $\tilde{u}(t,\theta)$ have the same type of zeros. The chain rule implies
\begin{equation*}
    \tilde{u_t}=\tilde{u}_{\theta\theta}+\frac{\tilde{u}_{\theta}}{\theta}+ \tilde{c}(t,\theta)\tilde{u}
\end{equation*}
for $\theta\in[0,\theta_2]$, where $\tilde{c}(t,\theta):=c(t,\theta)-\frac{b(\theta)}{2\theta}+\frac{b^2(\theta)}{4}-\frac{b_\theta(\theta)}{2}$. Note the term $\frac{b(\theta)}{\theta}$ is not singular at $\theta=0$ due to the nature of $b(\theta)$, that is, its first order term is $b_0\theta$.

Next, the rescaled solution will be cut off. Define the cut off function $\eta:[0,\pi]\rightarrow [0,1]$ given by
\begin{align*}
    \eta(\theta):=
	\begin{cases}
		 1 & \text{ for } \theta \in[0,\theta_1]\\
		 0<\eta(\theta) <1 & \text{ for } \theta \in(\theta_1,\theta_2)\\
		 0 & \text{ for } \theta \in[\theta_2,\pi]
 	\end{cases}
\end{align*}
which transitions smoothly from 1 to 0. 

Let $v:[t_1,t_2]\times[0,\pi]\rightarrow \mathbb{R}$ be defined by
\begin{align*}
    v(t,\theta):=
	\begin{cases}
		 \eta(\theta)[\tilde{u}(t,\theta)-1]+1 & \text{ for } \theta \in[0,\theta_2]\\
		 1 & \text{ for } \theta \in(\theta_2,\pi].
 	\end{cases}
\end{align*}

That is, $v(t,\theta)=\tilde{u}(t,\theta)$ for $\theta\in[0,\theta_1]$. For $\theta\in[\theta_1,\theta_2]$ there is a transition phase from $\tilde{u}$ to the constant function $1$. For $\theta\in[\theta_2,\pi]$, the singularity at $\theta=\pi$ does not play a role anymore, since $v(t,\theta)\equiv 1$ satisfies a trivial equation.

The chain rule says that $v(t,\theta)$ satisfies
\begin{equation*}
    v_t=v_{\theta\theta}+\frac{v_{\theta}}{\theta}-\frac{\eta_\theta[\tilde{u}-1]}{\theta}+\eta \tilde{c}\tilde{u}-\eta_{\theta\theta}[\tilde{u}-1]-2\eta_\theta\tilde{u}_\theta .
\end{equation*}

Now $d(t,\theta)$ is defined so that $v(t,\theta)$ satisfies the desired equation \eqref{cutoffPDE}. For $\theta \in[0,\theta_1]$, the only term that does not vanish is $\tilde{c}\tilde{u}$, since $\eta\equiv 1$ and $\eta_\theta\equiv0\equiv\eta_{\theta\theta}$. This defines $d$ in this interval. For $\theta \in(\theta_1,\theta_2)$, define most terms on the right hand side by $d(t,\theta)v$, as below. For $\theta \in[\theta_2,\pi]$, the function $v\equiv 1$ and $v_t=v_{\theta\theta}=\frac{v_{\theta}}{\theta}=0$. Hence, it satisfies a trivial equation and define $d:=0$. More precisely,
\begin{align*}
    d(t,\theta):=
	\begin{cases}
		 \tilde{c} & \text{ for } \theta \in[0,\theta_1]\\
		 \frac{1}{v}[-\frac{\eta_\theta[\tilde{u}-1]}{\theta}+\eta \tilde{c}\tilde{u}-\eta_{\theta\theta}[\tilde{u}-1]-2\eta_\theta\tilde{u}_\theta] & \text{ for } \theta \in(\theta_1,\theta_2)\\
		 0 &\text{ for } \theta \in[\theta_2,\pi]
 	\end{cases}
\end{align*}
is bounded, since all terms $\tilde{u},\tilde{u}_\theta,\eta,\eta_\theta,\eta_{\theta\theta},\tilde{c}$ are bounded for $\theta \in[0,\theta_2]$. Also, note $v>0$ for $\theta \in[\theta_1,\theta_2]$ and hence $\frac{1}{v}$ is well defined and bounded. Indeed, the solution $u$ is positive in this interval, and so is $\tilde{u}$, since they have the same type of zeros. If $\tilde{u}\geq 1$ it is clear that $v>0$ by its definition, and if $1>\tilde{u}>0$, one also obtains that $v>0$ by noticing that $\eta\in [0,1]$ for $\theta \in[0,\theta_2]$.

Hence, we have defined $v$ and $d$ satisfying \eqref{cutoffPDE} such that $v$ and $u$ have the same type of zeros and $v\equiv 1$ for $\theta \in[\theta_2,\pi]$.
\begin{flushright}
	$\blacksquare$
\end{flushright}
\end{pf}

In order to apply the dropping lemma to functions $v(t,\theta)$ satisfying the equation \eqref{cutoffPDE}, as in \cite{ChenPolacik96}, one still needs two adaptations. Firstly, the dropping lemma is proved for $\theta\in[0,1]$ and this can be circumvented by stretching the interval through $\theta\mapsto\pi\theta$. Secondly, in \cite{ChenPolacik96} it is considered Dirichlet boundary condition at the regular boundary $\theta=1$, but their proof works similarly for the Neumann case by changing the odd reflection done at the regular boundary $\theta=1$ to an even reflection. Such choice of reflections is done explicitly in \cite{Angenent88}, for different boundary conditions.

\begin{pf}\textbf{of Lemma \ref{droplem} (dropping lemma)} Firstly, we prove that $u$ has finitely many zeros. The Lemma \ref{adaptlem} implies that one can construct a $v$ satisfying \eqref{cutoffPDE} with same type of zeros of $u$. Due to the dropping Lemma in \cite{ChenPolacik96}, $v$ has finitely many zeros and consequently $u$ has finitely many zeros for $\theta\in[0,\theta_1]$. 

To conclude that $u$ also has finitely many zeros for $\theta\in[\theta_1,\pi]$, consider the change of coordinates $\tilde{\theta}:= \pi-\theta$. The solution $u(t,\tilde{\theta})$ satisfies the equation \eqref{cutoffPDE} with $\tilde{\theta}\in[0,\pi-\theta_1]$, and by the dropping lemma in \cite{ChenPolacik96}, it also has finitely many zeros for $\tilde{\theta}\in[0,\pi-\theta_1]$. Equivalently, $u$ has finitely many zeros for ${\theta}\in[\theta_1,\pi]$.

Secondly, we prove that multiple zeros must drop. Suppose $(t_0,\theta_0)$ is a multiple zero of a solution $u\not\equiv 0$ of \eqref{linaxisymPDE}. By the Lemma \ref{adaptlem}, there is a function $v(t,\theta)$ having zeros of the same type as $u(t,\theta)$ for $\theta\in[0,\theta_1]$ and some fixed $\theta_1\in(0,\pi)$. 

If $\theta_0\leq\theta_1$, then the dropping lemma in \cite{ChenPolacik96} implies that the number of zeros of $v(t,\theta)$ should drop. Since $v(t,\theta)$ is not zero for $\theta\in[\theta_1,\pi]$, then the zero that dropped should have occured for $\theta\in[0,\theta_1]$. This implies that some zero of $u(t,\theta)$ must have dropped, since they have the same type of zeros. 

If $\theta_0>\theta_1$, then consider the change of coordinates $\tilde{\theta}:=\pi-\theta$ and the same arguments as above show that the multiple zero of $u(t,\tilde{\theta})$ must have dropped for $\tilde{\theta}\in[0,\pi-\theta_1]$.

Thirdly, we prove that the zero number is not increasing in time. We already know that it must drop at multiple zeros. Suppose $(t_0,\theta_0)$ is a simple zero, that is $u(t_0,\theta_0)=0$ and $u_\theta(t_0,\theta_0)\neq 0$. Hence, the implicit function theorem says that $u(t,\theta(t))=0$ for an unique curve $\theta(t)$ in small neighborhood of $t_0$ such that $\theta(t_0)=\theta_0$. Hence, the simple zero persists and no new zeros are created.

\begin{flushright}
	$\blacksquare$
\end{flushright}
\end{pf}


\subsubsection{Proof 2}

This proof is an adaptation of Angenent \cite{Angenent88}, by rescaling the solution nearby a multiple zero of multiplicity $n$ and showing that there are $n$ zero curves backwards in time, and less curves forwards in time. We give a sketch of the proof.

For $t_0>0$, the \emph{localization} of the solution ${v}(t,\theta)$ of \eqref{linaxisymPDE} nearby the multiple zero $(t_0,\theta_0)$,
\begin{equation*}
    {w}(\tau,\xi):=e^{-\frac{\xi^2}{2}}{v}(t_0-e^{-2\tau},\theta_0+2e^{-\tau}\xi)
\end{equation*}
for $\tau\geq -\frac{1}{2}\log(t_0)=:\tau_0$. Due to the properly chosen parabolic rescaling, ${w}(\tau,\xi)$ satisfies
\begin{equation*}
    {w}_\tau=\frac{1}{2}{w}_{\xi\xi}+\frac{1}{2\tan (\theta_0+2e^{-\tau}\xi)}{w}_{\xi}-\frac{1}{2}(\xi^2-1){w}+q(\tau,\xi){w}
\end{equation*}
where $(\tau,\xi)\in(\tau_0,\infty)\times \mathbb{R}$ and $q(\tau,\xi)$ is bounded and decay with $\tau$. 

There are two cases: either the multiple zero is in the interior $\theta_0\in (0,\pi)$ or in one of the boundaries $\theta_0=0,\pi$. 

In the first case, the tangent term is regular and one can rescale this $w_\xi$ term out by an appropriate multiplying $w$ by an appropriate exponential. Then the arguments of Angenent \cite{Angenent88} hold. 

In the second case, there is a singular term only at one of the boundaries it is being zoomed in. One can reflect solutions along the other boundary, which is regular, and rescale the bounded terms to obtain 
\begin{equation*}
    {w}_\tau=\frac{1}{2}{w}_{\xi\xi}+\frac{1}{2\xi}{w}_{\xi}-\frac{1}{2}(\xi^2-1){w}+q(\tau,\xi){w}
\end{equation*}
for $x\in\mathbb{R}_+$. 

The operator $\frac{1}{2}{w}_{\xi\xi}+\frac{1}{2\xi}{w}_{\xi}$ is self-adjoint in $L^2_\xi([0,\infty))$ with weigth $\xi$. Due to Sturm-Liouville, the spectrum of such operator consists of simple eigenvalues and respective eigenfunctions $\phi_n(\xi)=e^{-\xi^2/2}L_n(\xi)$, where $L_n$ is a multiple of the $n$-th Laguerre polynomial. This eigenvalue problem is also known in the literature as the quantum harmonic oscillator in spherical coordinates. One can then follow the proof of Angenent by simply changing the functional spaces and its basis. 

\subsubsection{Consequences of the dropping lemma}

Two results follow by combining the dropping lemma \ref{droplem} and the asymptotic description in Theorem \ref{ConvEF}. The first is a result relating the zero number within invariant manifold and the Morse indices of equilibria. The second is the Morse-Smale property.
\begin{thm} \label{Znuminvmfld}
    \emph{\textbf{Zero number within Invariant Manifolds} \cite{Tabata80}, \cite{FiedlerBrunovsky86} }
    
    Consider a equilibria $u_\pm\in\mathcal{E}$ and a trajectory $u(t)\not\equiv u_\pm$ of \eqref{axisymPDE}. Then,
    \begin{enumerate}
    \item If $u(t) \in W^u(u_-)$, then $i(u_-)>z^t(u-u_-)$.
    \item If $u(t) \in W^s_{loc}(u_+)$, then $z^t(u-u_+)\geq i(u_+)$.
    \item If $u(t) \in  W^u(u_-)\cap W^s_{loc}(u_+)$, then            \begin{equation*}
            i(u_+)\leq z^t(u-u_\pm)< i(u_-). 
        \end{equation*}   
    \end{enumerate}
    These results also hold by replacing $u(t)-u_*$ with the tangent vector $u_t$.
\end{thm}

The above theorem implies that \eqref{axisymPDE} has no homoclinic orbits. Indeed, if there were any, then $i(u_*)<i(u_*)$, which is a contradiction.

    

This last theorem implies that if the semigroup has a finite number of equilibria, in which all are hyperbolic, then it is a Morse-Smale system in the sense of \cite{HaleMagalhaesOliva84}. Note that this property can hold even in case the equilibria are not hyperbolic, as in \cite{Henry85}.

\subsection{Sturm global structure}\label{sec:globalsturm}

This section gathers all the tools developed in the previous sections in order to construct the attractor for the parabolic equation with singular coefficients \eqref{axisymPDE} and prove the second part of the main Theorem \ref{attractorthm}.

Its proof is a consequence of two propositions. Firstly, due to the \emph{cascading principle}, it is enough to construct all heteroclinics between equilibria such that their Morse indices differ by 1. Secondly, on one direction, the \emph{blocking principle}: some conditions imply that there does not exist a heteroclinic connection; on the other direction, the \emph{liberalism principle}: if those conditions are violated, then there exists a heteroclinic.

The cascading and blocking principles follow from the dropping lemma and Morse-Smale property from Section \ref{sec:drop}, and we give a sketch as in \cite{FiedlerRocha96}. There is only a mild modification in the proof of the liberalism principle in Proposition \ref{estabhets}.

\begin{prop}\emph{\textbf{Cascading Principle} \cite{FiedlerRocha96}}\label{cascading}
 
 There exists a heteroclinic between two equilibria $u_\pm$ such that  $n:=i(u_-)-i(u_+)>0$ if, and only if, there exists a sequence (cascade) of equilibria $\{ v_k\}_{k=0}^n$ with $v_0:=u_-$ and $v_n:=u_+$, such that the following holds for all $k=0,...,n-1$
 \begin{enumerate}
     \item $i(v_{k+1})=i(v_k)+1$
     \item There exists a heteroclinic from $v_{k+1}$ to $v_k$
 \end{enumerate} 
\end{prop}

\begin{prop} \emph{\textbf{Blocking and Liberalism Principles} \cite{FiedlerRocha96}} \label{estabhets}

There exists a heteroclinic between the equilibria $v_{k+1}$ and $v_k$ with $i(v_{k+1})=i(v_k)+1$ if, and only if,
  \begin{enumerate}
     \item \emph{Morse permit: } $z(v_{k+1}-v_k)=i(v_k)$,
     \item \emph{Zero number permit: } $z(v_{k+1}-u_*)\neq z(v_k-u_*)$ for all equilibria $u_*$ between $v_{k+1}$ and $v_k$ along $M^u_\theta$ for some $\theta\in [0,\pi]$.
 \end{enumerate} 
\end{prop}

The blocking and liberalism principles assert that the Morse indices $i(.)$ and zero numbers $z(.)$ construct the global structure of the attractor explicitly. Those numbers can be obtained through the Sturm permutation, as in Section \ref{sec:perm}. 

In particular, one can check the zero number blocking for $\theta=0$ as it is done in \cite{FiedlerRocha96}. We prefer to state the condition for some $\theta\in [0,\pi]$ because the Sturm permutation in Section \ref{sec:perm} labels the equilibria along $M^u_\theta$ and $M^s_\theta$ for some $\theta\in [0,\pi]$. Moreover, those curves are computed for $\theta=\pi/2$ for the Chafee-Infante example in Section \ref{sec:CI} 

We now show that $u_*$ lies in between $u_-$ and $u_+$ at $\theta=0$ if, and only if it is also between $u_\pm$ along $M^u_\theta$ for any $\theta\in [0,\pi]$. Indeed, due to continuity with respect to the initial data $(0,a,0)\in\mathbb{R}^3$ of the shooting flow \eqref{shootflow}, the curve $M^u_{\theta}$ for fixed $\theta\in [0,\pi)$ is continuous and the order of $a\in\mathbb{R}$ induces an order along $M^u_\theta$, hence the parametrization respects its labeling. At $\theta=\pi$, continuity also yields an ordering of the equilibria within $M^u_{\theta}$.

Note one can replace $M^u_\theta$ in the zero number blocking by $M^s_\theta$, since similar arguments as above hold and show that $u_*$ lies in between $u_-$ and $u_+$ at $\theta=\pi$ if, and only if it is also between $u_\pm$ along $M^s_\theta$ for some $\theta\in [0,\pi]$. 

The two propositions above yield the existence of heteroclinics between $u_-$ and $u_+$ if they are \emph{cascadly adjacent}, namely, if there exists a cascade of equilibria $\{ v_k\}_{k=0}^n$ with $v_0:=u_-$ and $v_n:=u_+$ such that for all $k=0,...,n-1$ the following conditions hold:
  \begin{enumerate}
     \item $i(v_{k+1})=i(v_k)+1$,
     \item $z(v_k-v_{k+1})=i(v_{k+1})$,
     \item $z(v_{k+1}-u_*)\neq z(v_k-u_*)$ for all equilibria $u_*$ between $v_{k+1}$ and $v_k$ along $M^u_\theta$ for some $\theta\in [0,\pi]$.
 \end{enumerate} 
 
On the other hand, the main Theorem \ref{attractorthm} yields a result through the notion of adjacency in the introduction, which does not involve a cascade. These notions of adjacency coincide, and this is the core of Wolfrum's ideas in \cite{Wolfrum02}. 

\begin{prop} \emph{\textbf{Wolfrum's equivalence}} \label{wolfrumlemma}
Consider two equilibria $u_\pm\in\mathcal{E}$ such that $n:=i(u_-)-i(u_+)>0$. The equilibria $u_\pm$ are adjacent if, and only if they are cascadly adjacent.
\end{prop}

The proof of the cascading proposition \ref{cascading} follows \cite{FiedlerRocha96} word by word.

For the proof of the liberalism theorem, it is used the Conley index to detect orbits between $u_-$ and $u_+$. We give a brief introduction of Conley's theory, and how it can be applied in this context. See Chapters 22 to 24 in \cite{Smoller83} for a brief account of the Conley index, and its extension to infinite dimensional systems in \cite{Rybakowski82}.

Consider the space $\mathcal{X}$ of all topological spaces and the equivalence relation given by $Y\sim Z$ for $Y,Z\in \mathcal{X}$ if, and only if $Y$ is homotopy equivalent to $Z$, that is, there are continuous maps $f:Y\rightarrow Z$ and $g:Z\rightarrow Y$ such that $f\circ g$ and $g\circ f$ are homotopic to $id_Z$ and $id_Y$, respectively. Then, the quotient space $\mathcal{X}/ \sim$ describes the homotopy equivalent classes $[Y]$ of all topological spaces which have the same homotopy type. Intuitively, $[Y]$ describes all topological spaces which can be continuously deformed into $Y$.

Suppose $\Sigma$ is an invariant isolated set, that is, it is invariant with respect to positive and negative time of the semiflow, and it has a closed neighborhood $N$ such that $\Sigma$ is contained in the interior of $N$ with $\Sigma$ being the maximal invariant subset of $N$. 

Denote $\partial_e N\subset \partial N$ the \emph{exit set of $N$}, that is, the set of points which are not strict ingressing in $N$,
\begin{equation*}
    \partial_e N:=\{u_0\in N \text{ $|$ } u(t)\not\in N \text{ for all sufficiently small $t>0$}\} .
\end{equation*}

The \emph{Conley index} is defined as
\begin{equation*}
    C(\Sigma):=[N/ \partial_e N]
\end{equation*}
namely the homotopy equivalent class of the quotient space of the isolating neighborhood $N$ relative to its exit set $\partial_e N$. Such index is homotopy invariant and does not depend on the particular choice of isolating neighborhood $N$.

We compute the Conley index for two examples.

Firstly, the Conley index of a hyperbolic equilibria $u_+$ with Morse index $n$. Consider a closed ball $N\subset X^\alpha$ centered at $u_+$ without any other equilibria in $N$, as isolating neighborhood. The flow provides a homotopy that contracts along the stable directions to the equilibria $u_+$. 
Then, $N$ is homotoped to a $n$-dimensional ball $B^n$ in the finite dimensional space spanned by the first $n$ eigenfunctions, related to the unstable directions. Note the exit set $\partial_e B^n=\partial B^n=\mathbb{S}^{n-1}$, since after the homotopy there is no more stable direction and the equilibria is hyperbolic. Therefore, the quotient of a $n$-ball and its boundary is an $n$-sphere,
\begin{equation*}
    C(u_+)=[N/\partial_e N]=[B^n/\partial_e B^n]=[B^n/\mathbb{S}^{n-1}]=[\mathbb{S}^n].
\end{equation*} 

Secondly, the Conley index of the union of two disjoint invariant sets, for example $u_-$ and $u_+$ with respective disjoint isolating neighborhoods $N_-$ and $N_+$. Then, $N_-\cup N_+$ is an isolating neighborhood of $\{ u_-,u_+\}$. By definition of the wedge sum 
\begin{align*}
    C(\{ u_-,u_+\})&=\left[\frac{N_-\cup N_+}{\partial_e (N_-\cup N_+)}\right]\\
    &=\left[\frac{N_-}{\partial_e N_-}\vee\frac{N_+}{\partial_e N_+}\right]=C(u_-)\vee C(u_+).
\end{align*}

The Conley index can be applied to detect heteroclinics as follows. Construct a closed neighborhood $N$ such that its maximal invariant subspace is the closure of the set of heteroclinics between $u_\pm$,
\begin{equation*}
    \Sigma=\{u_-,u_+\}\cup \overline{W^u(u_-)\cap W^s(u_+)}.
\end{equation*} 

Suppose, towards a contradiction, that there are no heteroclinics connecting $u_-$ and $u_+$, that is, $\Sigma=\{ u_-,u_+\}$. Then, the index is given by the wedge sum $C(\Sigma)=[\mathbb{S}^n]\vee [\mathbb{S}^m]$, where $n,m$ are the respective Morse index of $u_-$ and $u_+$.

If, on the other hand, one can prove that $C(\Sigma)=[0]$, where $[0]$ means that the index is given by the homotopy equivalent class of a point, this would yield a contradiction and there should be a connection between $u_-$ and $u_+$. Moreover, the Morse-Smale structure excludes connection from $u_+$ to $u_-$, and hence there is a connection from $u_-$ to $u_+$.

Hence, there are three ingredients missing in the proof: the Conley index can be applied at all, the construction of a isolating neighborhood $N$ of $\Sigma$ and the proof that $C(\Sigma)=[0]$.

\begin{pf}\textbf{of Proposition \ref{estabhets}}

($\implies$) This part is called blocking and has same proof as in \cite{FiedlerRocha96}.

($\impliedby$) This is also called liberalism in \cite{FiedlerRocha96}. Consider hyperbolic equilibria $u_-,u_+$ such that $i(u_-)=i(u_+)+1$ and satisfies both the Morse and the zero number permit conditions. Without loss of generality, assume $u_-(0)>u_+(0)$. 

It is used the Conley index to detect orbits between $u_-$ and $u_+$. Note that the semiflow generated by the equation \eqref{axisymPDE} on the Banach space $X^\alpha$ is admissible for the Conley index theory in the sense of \cite{Rybakowski82}, due to a compactness property that is satisfied by the parabolic equation \eqref{axisymPDE}, namely that trajectories are precompact in phase space. See Theorem 3.3.6 in \cite{Henry81}.

As mentioned above, in order to apply the Conley index concepts we need to construct appropriate neighborhoods and show that the Conley index is $[0]$. 

Consider the closed set
\begin{equation*}
    K({u_\pm}):= \left\{ u \in X^\alpha \mid 
        \begin{array}{c} 
        z(u-u_-)=i(u_+)=z(u-u_+) \\
        u_+(0)\leq u(0)\leq u_-(0) 
        \end{array} \right\}
\end{equation*}

Consider also closed $\epsilon$-balls $B_\epsilon(u_\pm)$ centered at $u_\pm$ such that they do not have any other equilibria besides $u_\pm$, respectively, for some $\epsilon>0$. 

Define
\begin{equation*}
    N_\epsilon (u_\pm):=B_\epsilon(u_-)\cup B_\epsilon(u_+) \cup K({u_\pm}).    
\end{equation*}

The zero number blocking condition implies there are no equilibria in $K({u_\pm})$ besides possibly $u_-$ and $u_+$. Hence, $N_\epsilon (u_\pm)$ also has no equilibria besides $u_-$ and $u_+$. 

Denote $\Sigma$ the maximal invariant subset of $N_\epsilon$. We claim that $\Sigma$ is the set of the heteroclinics from $u_-$ to $u_+$ given by $\overline{W^u(u_-)\cap W^s(u_-)}$, and the equilibria themselves. 

On one hand, since $\Sigma$ is globally invariant, then it is contained in the attractor $\mathcal{A}$, which consists of equilibria and heteroclinics. Since there are no other equilibria in $N_\epsilon (u_\pm)$ besides $u_\pm$, then the only heteroclinics that can occur are between them.

On the other hand, Theorem \ref{Znuminvmfld} implies that along a heteroclinic $u(t)\in\mathcal{H}$ the zero number satisfies $z^t(u-u_\pm)=i(u_+)$ for all time, since $i(u_-)=i(u_+)+1$. Therefore $u(t)\in K({u_\pm})$ and the closure of the orbit is contained in $N_\epsilon (u_\pm)$. Since the closure of the heteroclinic is invariant, it must be contained in $\Sigma$.

Lastly, it is proven that $C(\Sigma)=[0]$ in three steps, yielding the desired contradiction and the proof of the theorem. We modify the first and second step from \cite{FiedlerRocha96}, whereas the third remain the same. 

In the first step, a model is constructed displaying a saddle-node bifurcation with respect to a parameter $\mu$, for $n:=z(u_+-u_-)\in\mathbb{N}$ fixed,
\begin{equation}\label{prototype}
    v_t=a(\xi,v,v_\xi)\left[\left(v_{\xi\xi}+\frac{1}{\tan(\xi)}v_\xi\right) +\lambda_n v\right]+ g_n(\mu,\xi,v)
\end{equation}
where $\xi\in [0,\pi]$ has Neumann boundary conditions, $\lambda_n=n(n+1)$ are the eigenvalues of the axisymmetric laplacian with the Legendre polynomials $P_n(cos(\xi))$ as eigenfunctions, and 
\begin{equation*}
    g_n(\mu,\xi,v):=\left[v^2-\mu P_n^2\right]P_n.
\end{equation*}



For $\mu>0$, the equilibria solution of \eqref{prototype} are $v_\pm=\pm \sqrt{\mu}P_n(\cos(\xi))$, since $P_n$ are the eigenfunctions of the axially symmetric Laplacian. Furthermore, we have 
\begin{equation}
    z(v_+-v_-)=n
\end{equation}
since the $n$ intersections of $v_-$ and $v_+$ will be at its $n$ zeroes. 

Moreover, $v_\pm$ are hyperbolic equilibria for small $\mu>0$, such that $i(v_+)=n+1$ and $i(v_-)=n$. Indeed, parametrize the bifurcating branches by $\mu=s^2$ so that $v(s,\xi)=sP_n(\cos(\xi))$, where $s>0$ correspond to $v_+$ and $s<0$ to $v_-$. Linearizing at the equilibrium $v_\pm$ yields the following linear operator
\begin{equation*}
    T_n(s)v:=a(\xi,sP_n,-sP'_n\cdot \sin(\xi))\left[v_{\xi\xi}+\frac{1}{\tan(\xi)}v_\xi +\lambda_n v\right]+2svP^2_n.
\end{equation*}

This operator can be seen as a Sturm-Liouville eigenvalue problem in the space $L^2_w$ with appropriate weight $w:=P^2_n(\cos(\xi))$, namely, $T_n(s)v=\eta P^2_n v$. Notice that for $v=P_n(\cos(\xi))$ the first term vanish, and hence the eigenproblem becomes 
\begin{equation*}
    2sP^2_nv=\eta P_n^2v.
\end{equation*}


Hence, for each $n$ fixed, $\eta_n(s)=2s$ is an eigenvalue with $P_n(\cos(\xi))$ its corresponding eigenfunction, since the terms inside the brackets yield the eigenvalue problem for the axisymmetric Laplacian and vanish. 

We now use a perturbation argument in Sturm-Liouville theory. For $\mu=s^2=0$, the eigenvalues of $L_n(0)$ in $L^2_w$ coincide with the eigenvalues of the usual axisymmetric laplacian such that there is one eigenvalue $\eta_n(0)=0$ and $n$ positive eigenvalues. For small $\mu<0$, the number of positive eigenvalues persist, and there is no eigenvalue $0$, since $\eta(\mu)<0$; whereas for small $\mu>0$, the number of positive eigenvalues increases by 1, and there is no eigenvalue $0$, since $\eta(\mu)>0$. This yields the desired claim about hyperbolicity and the Morse index.

Now consider the semilinear parabolic equation such that \eqref{prototype} is its equilibria equation. The equilibria $v_\pm$ together with their connecting orbits of the corresponding evolution equation form an isolated invariant set
\begin{equation*}
    \Sigma_\mu(v_\pm):= \overline{W^u(v_-)\cap W^s(v_+)}
\end{equation*}
with isolating neighborhood $N_\epsilon(v_\pm)$, and the bifurcation parameter can also be seen as a homotopy parameter. Hence the Conley index is of a point by homotopy invariance as desired, that is,
\begin{equation}\label{conley1}
    C(\Sigma_\mu(v_\pm))=C(\Sigma_0(v_\pm))=[0].
\end{equation}

In the second step, the equilibria $v_-$ and $v_+$ are transformed respectively into $u_-$ and $u_+$ via a diffeomorphism which is not a homotopy. 

Recall $n=z(v_--v_+)=z(u_+-u_-)$. Hence, choose $\xi(\theta)$ a smooth diffeomorphism of $[0,\pi]$ that maps the zeros of $v_-(\xi)-v_+(\xi)$ to the zeros of $u_-(\theta)-u_+(\theta)$. Therefore, from now on we suppose that the zeros of $v_-(\xi(\theta))-v_+(\xi(\theta))$ and $u_-(\theta)-u_+(\theta)$ occur in the same points in $\theta\in [0,\pi]$. From now on, we write the unknown $v$ as $v(\theta)$ when we actually mean $v(\xi(\theta))$, in order to simplify the notation.

Consider the transformation 
\begin{align*}
    L: X^\alpha &\to X^\alpha   \\ 
    v&\mapsto l(\theta)[v(\theta)-v_-(\theta)]+u_-(\theta)
\end{align*}
where $l(\theta)$ is defined pointwise through
\begin{align*}
    l(\theta):=
    \begin{cases}
        \frac{u_+(\theta)-u_-(\theta)}{v_+(\theta)-v_-(\theta)} &, \text{ if } v_+(\theta)\neq v_-(\theta) \\
        \frac{\partial_\theta(u_{+}(\theta)-u_{-}(\theta))}{\partial_\theta(v_{+}(\theta)-v_{-}(\theta))} &, \text{ if } v_+(\theta)= v_-(\theta)
    \end{cases}
\end{align*}
such that the coefficient $\alpha$ is smooth and nonzero due to the l'H\^opital rule. Hence, $L(v_-)=u_-$ and $L(v_+)=u_+$ as desired. Note we supposed $2\alpha+\beta>1$ so that solutions $u_\pm\in C^1$, hence $L$ is of this regularity as well. Moreover, $L$ is invertible with inverse having the same regularity. In particular, it is a homeomorphism, and hence a homotopy equivalence.

Moreover, the number of intersections of functions is invariant under the map $L$,
\begin{equation}
    z(L (v(\xi)-\tilde{v}(\xi)))=z(v(\theta)-\tilde{v}(\theta))
\end{equation}
and hence $K({v_\pm})$ is mapped to $K({u_\pm})$ under $L$. 

Consider $w(t,\theta):=L(v(t,\xi))$, hence the map $L$ modifies the equation \eqref{prototype} into the following equation
\begin{equation}\label{IDK}
    w_t=\tilde{a}(\theta,w,w_\theta)w_{\theta\theta}+\tilde{b}(\theta,w,w_\theta)\frac{w_\theta}{\tan(\theta)}+\tilde{f}(\theta,w,w_\theta)
\end{equation}
where the Neumann boundary conditions are preserved, and the terms $\tilde{a},\tilde{b},\tilde{f}$ are 
\begin{align*}
    \tilde{a}(\theta,w,w_\theta):=&\frac{\theta_\xi^2}{l(\theta)} \cdot a(\theta,L^{-1}(w),\partial_\theta L^{-1}(w))\\
    \tilde{b}(\theta,w,w_\theta):=&\frac{\theta_\xi}{l(\theta)} \cdot a(\theta,L^{-1}(w),\partial_\theta L^{-1}(w))\\    
    \tilde{f}(\theta,w,w_\theta):=&g_n(\mu,\theta,L^{-1}(w),\partial_\theta L^{-1}(w))+ (w_\theta\cdot \theta_{\xi\xi}-\partial_\xi^2u_-)-\frac{l_{\xi\xi} \cdot (w-u_-)_\xi}{l}\\
    &-\frac{l\partial_\xi u_- +l_{\xi} \cdot (w-u_-)}{l^2}+\partial_\xi v_--\lambda a(\theta,L^{-1}(w),\partial_\theta L^{-1}(w))\cdot L^{-1}(w).
\end{align*}

Note that the equilibria $v_\pm$ are mapped into $w_\pm:=L(v_\pm)=u_\pm$, which are equilibria of \eqref{IDK}, with same zero numbers and Morse indices as $v_\pm$ and $u_\pm$.

The isolated invariant set $\Sigma_\mu(v_\pm)$ is transformed into $L(\Sigma_\mu(v_\pm))=\Sigma_\mu(w_\pm)$, which is still isolated and invariant, with invariant neighborhood $L(N_\epsilon(v_\pm))=N_\epsilon(w_\pm)$. Moreover, the Conley index is preserved, since $L$ is a homotopy equivalence,
\begin{equation}\label{conley2}
   C(\Sigma_\mu(v_\pm))= C(L(\Sigma_\mu(v_\pm)))=C(\Sigma_\mu(w_\pm)).
\end{equation}

Hence, one identifies the equilibria $v_\pm$ in the model constructed \eqref{prototype} with the equilibria $w_\pm=u_\pm$ from the equation \eqref{axisymPDE}, by preserving neighborhoods and the Conley index, since $L$ is a homotopy equivalence. The identified equilibria $u_\pm$ satisfy the equation \eqref{IDK}, and we still have to modify it to become \eqref{axisymPDE}. For such, we perform now a last homotopy between the solutions $w$ and $u$.

In the third step, we homotope the diffusion coefficient $\tilde{a},\tilde{b}$ and nonlinearity $\tilde{f}$ from the equation \eqref{IDK} to be the desired diffusion $a$ and reaction $f$ from the equation \eqref{axisymPDE}. Indeed, consider the parabolic equation 
\begin{equation*}
    u_t=a^\tau(\theta,u,u_\theta)u_{\theta\theta}+b^\tau(\theta,u,u_\theta)\frac{u_\theta}{\tan(\theta)}+f^\tau(\theta,u,u_\theta)
\end{equation*}
where
\begin{align*}
    a^\tau&:=\tau \tilde{a}+(1-\tau)a+\sum_{i=- \text{ , } +}\chi_{u_i}\mu_{u_i}(\tau)[u-{u_i}(x)]\\
    b^\tau(\theta)&:=\tau \tilde{b}+(1-\tau)a+\sum_{i=- \text{ , } +}\chi_{u_i}\mu_{u_i}(\tau)[u-{u_i}(x)]\\
    f^\tau(x,u,u_x)&:=\tau \tilde{f}+(1-\tau)f+\sum_{i=- \text{ , } +}\chi_{u_i}\mu_{u_i}(\tau)[u-{u_i}(x)]
\end{align*}
and $\chi_{u_i}$ are cut-offs being 1 nearby $u_i$ and zero far away, the coefficients $\mu_i(\tau)$ are zero near $\tau=0$ and $1$ and shift the spectra of the linearization at $u_\pm$ such that uniform hyperbolicity of these equilibria is guaranteed during the homotopy. Note that $u_\pm$ have the same Morse indices, as solutions of both equations \eqref{axisymPDE} and \eqref{IDK}. Therefore, the $\mu_i(\tau)$ only makes sure none of these eigenvalues cross the imaginary axis.

Consider $u_\pm$ and their connecting orbits during this homotopy, 
\begin{equation*}
    \Sigma^\tau(u_\pm):=\overline{W^u(u_-)\cap W^u(u_+)}.
\end{equation*}

Note that $\Sigma^\tau(u_\pm)\subseteq K({u_\pm})$, for all $\tau\in [0,1]$, since the dropping lemma holds throughout the homotopy. The equilibria $u_\pm$ do not bifurcate as $\tau$ changes, due to hyperbolicity. Choosing $\epsilon>0$ small enough, the neighborhoods $N_\epsilon(u_\pm)$ form an isolating neighborhood of $\Sigma^\tau(u_\pm)$ throughout the homotopy. Indeed, $\Sigma^\tau(u_\pm)$ can never touch the boundary of $K({u_\pm})$, except at the points $u_\pm$ by the dropping lemma. Once again the Conley index is preserved by homotopy invariance,
\begin{equation}\label{conley3}
    C(\Sigma(u_\pm))=C(\Sigma^0(u_\pm))=C(\Sigma^\tau(u_\pm))=C(\Sigma^1(u_\pm))=C(\Sigma_\mu(w_\pm)).
\end{equation}

Finally, the equations \eqref{conley1}, \eqref{conley2} and \eqref{conley3} yield that the Conley index of $\Sigma$ is the homotopy type of a point, and hence the desired result:
\begin{equation}\label{conley4}
    C(\Sigma(u_\pm))=C(\Sigma_\mu(w_\pm))=C(\Sigma_\mu(v_\pm))=[0].
\end{equation}
    
\begin{flushright}
	$\blacksquare$
\end{flushright}
\end{pf}

\section{Example: Axisymmetric Chafee-Infante}\label{sec:CI}

In this section it is given an example of the theory above, namely, it is constructed the attractor of the axially symmetric Chafee-Infante problem, 
\begin{equation}\label{PDECI}
    u_t= a(\theta,u,u_\theta)\left[ u_{\theta\theta}+\frac{1}{\tan(\theta)}u_\theta+\lambda u[1-u^2]\right]
\end{equation}
where $\theta\in[0,\pi]$ has Neumann boundary conditions, $a>0$, and initial data $u_0\in X^\alpha$, so that the equation generates a dynamical system in such space.

We will prove that the Sturm permutation for the axially symmetric case is the same as the permutation for the regular case in \cite{FiedlerRocha96}. Hence, both attractors are connection-wise the same. The only difference lies in the shape of equilibria, and that the parameter $\lambda$ must lie between two eigenvalues of the appropriate diffusion operator.

The proof is divided in the upcoming subsections. Following the shooting arguments in Section \ref{sec:perm}, we explicitly construct the shooting manifolds. Then we count how many times they intersect, yielding all equilibria, and if such intersections are transverse, yielding hyperbolicity. Lastly those intersection points are labeled accordingly, yielding the Sturm permutation $\sigma$, and hence the attractor $\mathcal{A}$.

The equilibria equation describing the shooting curves is
\begin{align}\label{shootCI}
\begin{cases}
    u_{\tau}&= p\\
	p_{\tau}&=-\lambda u[1-u^2]\sin^2(\theta)\\
	\theta_\tau&=\sin(\theta).
\end{cases}
\end{align}

Note that solutions of the shooting \eqref{shootCI} exist for all $\theta\in [0,\pi]$ and all initial data $d\in L_0:=\{ (0,d,0)\text{ $|$ } d\in\mathbb{R}\}$ or $e\in L_\pi:=\{ (\pi,e,0)\text{ $|$ } e\in\mathbb{R}\}$. Indeed, solutions of \eqref{shootCI} are bounded, since we can compare solutions of the axially symmetric Chafee-Infante \eqref{shootCI} with the usual Chafee-Infante, which is known to have global bounded solutions. 

This system possesses two symmetries, namely invariance under
\begin{align}
    \text{\emph{time reversal: }}& \tau \mapsto -\tau,\label{sym1} \\
    \text{\emph{reflection: }} & u\mapsto -u,\label{sym2}
\end{align}
where both symmetries also changes the sign of $p:=u_\tau$.
\subsection{Construction of the shooting curves}

The stable shooting manifold $M^s$ is obtained through the time reversal \eqref{sym1}, which is simply a reflection in the $p$-axis of the unstable shooting manifold $M^u$.

In order to construct the unstable shooting manifold $M^u$, we analyze four regions for the initial data $(0,d,0)\in\mathbb{R}^3$ constrained to the trivial equilibria $d\equiv-1,0,1$, for all parameter values $\lambda>0$. 

Note that part of the unstable shooting manifold $M^u$, namely when $d<0$, is obtained through a rotation by $\pi$, fixing the origin, of the piece of the shooting manifold $M^u$ when $d>0$, due to the reflectional symmetry \eqref{sym2}. 

If $d>1$, then the corresponding solution remains bigger than 1 for small time by continuity. Hence, the shooting flow \eqref{shootCI} implies that $p_\tau>0$ and the shooting manifold $M^u|_{d>1}$ increases in the $p$ direction as $\theta$ increases. 

For $d\in (0,1)$, we will show that the unstable shooting manifold $M^u|_{d\in(0,1)}$ winds around the trivial equilibrium $d\equiv 0$. More precisely, the angle and radius of the shooting manifold in polar coordinates are monotone with respect to its parametrization given by the initial data $d\in\mathbb{R}$.

This was proved in \cite{RochaHale85} using the Hamiltonian structure of the Chafee-Infante system, which can not be applied for the system \eqref{shootCI}, since it is nonautonomous. Instead, we adapt ideas of \cite{Lebovitz}. 

Indeed, the shooting flow \eqref{shootCI} in polar coordinates with the clockwise angle, $(u,p)=:(\rho\cos(\mu),-\rho\sin(\mu))$, is given by
\begin{align}\label{polarCI}
\begin{cases}
    \rho_\tau&=\rho \sin(\mu)\cos(\mu)\left[ \lambda (1-{\rho}^2 \cos^2(\mu)) \sin^2(\theta)-1\right]\\
    \mu_\tau&= \sin^2(\mu)+ \lambda [1-{\rho}^2 \cos^2(\mu)] \sin^2(\theta) \cos^2(\mu)\\
	\theta_\tau&=\sin(\theta)
\end{cases}
\end{align}
with $\displaystyle \lim_{\tau=-\infty} \mu(\lambda,\tau)=0$ describing $L_0$. Note that 
\begin{equation}\label{taumonot}
    \mu_\tau>0    
\end{equation} 
for $|\rho|<1$, that is, the angle $\mu$ is increasing in $\tau$ and each solution within the shooting manifolds are winding clockwise around the trivial equilibria $0$ as $\tau$ increases.

Consider the Lipschitz map $F(\lambda,\theta,\rho,\mu):R^4\to R^3$, where each coordinate $F_i$ correspond to the $i$-th line of the right-hand side in \eqref{polarCI}. 

We show the monotonicity of the angle $\mu$ with respect to the initial data $d\in(0,1)$, that is, the angle $\mu$ decreases as $d$ increases. This means that the bigger the initial data $d\in(0,1)$, smaller the angle, hence outer orbits rotate slower than inner orbits. 

\begin{lem}
    Let $(\rho,\mu)$ and $(\tilde{\rho},\tilde{\mu})$ be solutions of \eqref{polarCI} with different initial data given by $\lim_{\tau\to-\infty}(\rho(\lambda,\tau),\mu(\lambda,\tau))=(d,0)$ and $\lim_{\tau\to-\infty}(\tilde{\rho}(\lambda,\tau),\tilde{\mu}(\lambda,\tau))=(\tilde{d},0)$ with ordering $0<d<\tilde{d}<1$. Then
    \begin{equation}\label{mumonot}
        \mu(\lambda,\tau)>\tilde{\mu}(\lambda,\tau)    
    \end{equation}
    and
        \begin{equation}\label{rhomonot}
        \rho(\lambda,\tau)< \tilde{\rho}(\lambda,\tau)    
    \end{equation}
    for all $\tau\in\mathbb{R}$. 
    Moreover, if $\lambda>\tilde{\lambda}$ in \eqref{polarCI}, then
    \begin{equation}\label{lambdamonot}
        \mu(\lambda,\tau)>\mu(\tilde{\lambda},\tau)
    \end{equation}
    for all $\tau\in\mathbb{R}$ and fixed initial data $d\in (0,1)$. 
\end{lem}

From now on, we abuse the notation and suppress one of the coordinates of $\mu(\lambda,\tau)$, remarking only the one of importance in such equation, namely either $\mu(\lambda)$ or $\mu(\tau)$, even though $\mu$ still depends on those two variables. 
\begin{pf}
Firstly we show a weaker version of \eqref{mumonot} with a non strict inequality, namely
\begin{equation}\label{nonstrct}
    \mu(\tau)\geq\tilde{\mu}(\tau)    
\end{equation} 
for all $\tau\in(-\infty,\infty)$. 

Suppose, towards a contradiction, that 
\begin{equation}\label{contrtau1}
    \mu(\tau_1)<\tilde{\mu}(\tau_1)    
\end{equation} 
for some $\tau_1\in(-\infty,\infty)$. 

We show this inequality \eqref{contrtau1} also holds for $\tau\in(\tau_2,\tau_1]$ for some $\tau_2<\tau_1$. Note that for $\tau$ large and negative, the flow of the angle in \eqref{polarCI} is given by its linearization,
\begin{equation}\label{linshootCI}
    \mu_\tau= \sin^2(\mu)+ \lambda f_u(d,0) \sin^2(\theta) \cos^2(\mu)
\end{equation}
where $f_u(d,0)=1-3d^2$. The angle $\tilde{\mu}$ satisfies a similar equation with linearization given by $f_u(\tilde{d},0)=1-3\tilde{d}^2$. 

Indeed, nearby a non-hyperbolic fixed point, the flow \eqref{shootCI} is topologically equivalent to a decoupled system as in \cite{Shoshitaishvili71}, where the first equation describes the flow on the center manifold, and the second describes the linear hyperbolic dynamics. If the equilibria is hyperbolic, there is no center manifold and this breaks down to the Hartman-Grobman theorem. Since the shooting manifolds are the strong unstable and stable manifolds, there is no center direction within them, and the flow is topological equivalent to its corresponding hyperbolic part of the linearization. Therefore, we linearize \eqref{shootCI} at $(0,d,0)\in L_0$, then change to polar coordinates, yielding \eqref{linshootCI}.

Note that $f_u(d,0)>f_u(\tilde{d},0)$, since $0<d<\tilde{d}<1$. By the comparison theorem in \cite{Bundincevic10}, one obtains that for such linearizations, 
\begin{equation}\label{lintaustar}
    \mu(\tau)>\tilde{\mu}(\tau)    
\end{equation}
for all $\tau\in(-\infty,\tau^*)$ with $\tau^*$ negative and large such that $\tau^*<\tau_1$, that is, so that the nonlinear system \eqref{polarCI} for $\mu$ is topological equivalent to the linear one \eqref{linshootCI}.

By the intermediate value theorem, there exists $\tau_2\in(-\tau^*,\tau_1)$ such that $\mu(\tau_2)=\tilde{\mu}(\tau_2)$. We choose the biggest of those values, due to continuity of $\mu$ up to $\tau_1$, yielding
\begin{equation}\label{monot2}
    \mu(\tau)<\tilde{\mu}(\tau)    
\end{equation} 
for $\tau\in(\tau_2,\tau_1]$, which extends the inequality \eqref{contrtau1} as claimed.

On the other hand, the integral formulation of \eqref{polarCI} yields that
\begin{equation}\label{CUBA}
   \mu(\textcolor{black}{\tau})-\mu(\tau_2)= \int_{\tau_2}^{\textcolor{black}{\tau}} F_2(\lambda,\theta,\rho,\mu) d\textcolor{black}{s} 
\end{equation} 
\textcolor{black}{for any $\tau\in(\tau_2,\tau_1]$}. Similarly for $\tilde{\mu}$. 

Consider the difference $\tilde{\mu}-\mu$ of the above representation. Notice that $\mu(\tau_2)=\tilde{\mu}(\tau_2)$ and $F_2$ is Lispchitz in $\mu$ and $\rho$, while $\lambda$ and $\theta$ are fixed,
\begin{align*}
   |\tilde{\mu}(\textcolor{black}{\tau})-\mu(\textcolor{black}{\tau})|\leq c(\lambda,\theta) \int_{\tau_2}^{\textcolor{black}{\tau}} \sqrt{|\tilde{\rho}-\rho|^2+|\tilde{\mu}-\mu|^2} d\textcolor{black}{s} .
\end{align*}

Note that the square root of a sum is less than the sum of the square roots. 
Moreover, the solutions $\rho,\tilde{\rho}$ of \eqref{shootCI} are bounded, hence $|\tilde{\rho}-\rho|$ is bounded, say by $M$. Lastly, one can get rid of the norms in $|\tilde{\mu}-\mu|$, due to \eqref{monot2}. These considerations yield
\begin{align}\label{rhoormu}
   \tilde{\mu}(\textcolor{black}{\tau})-\mu(\textcolor{black}{\tau})\leq c\cdot M\int_{\tau_2}^{\textcolor{black}{\tau}} d\textcolor{black}{s} +c\int_{\tau_2}^{\textcolor{black}{\tau}}(\tilde{\mu}-\mu) d\textcolor{black}{s} .
\end{align}

The mean value theorem \textcolor{black}{for definite integrals guarantees there is a} $\tau_3\in(\tau_2,\textcolor{black}{\tau})$ so that
$\textcolor{black}{\tau}-\tau_2=\int_{\tau_2}^{\tau} (\tilde{\mu}-\mu)d\textcolor{black}{s}/(\tilde{\mu}(\tau_3)-\mu(\tau_3))$, 
where the denominator is well defined due to \eqref{monot2}. 
\textcolor{black}{Moreover, for $\epsilon>0$ sufficiently small, we let $m_\epsilon:=\inf_{s\in[\tau_2+\epsilon,\tau_1]} \tilde{\mu}(s)-\mu(s)$, which is well defined and bounded since $\tilde{\mu},\mu$ are continuous and the interval is compact. Also, $m_\epsilon>0$ due to the definition of $\tau_2$ and \eqref{monot2}. We plug in $\tau-\tau_2$ from such mean value formula together with the bound $1/(\tilde{\mu}(\tau_3)-\mu(\tau_3))\leq 1/m_\epsilon$ into \eqref{rhoormu}, yielding }
\begin{equation}\label{CUBA2}
   \tilde{\mu}(\textcolor{black}{\tau})-\mu(\textcolor{black}{\tau})\leq \left[\frac{c\cdot M}{m_\epsilon}+c\right]\int_{\tau_2}^{\textcolor{black}{\tau}} (\tilde{\mu}-\mu) d\textcolor{black}{s}
\end{equation}
\textcolor{black}{for any $\tau\in [\tau_2+\epsilon,\tau_1]$.}

The integral Gr\"{o}nwall inequality implies that $\tilde{\mu}(\tau)-\mu(\tau)\leq 0$ \textcolor{black}{for any $\tau\in [\tau_2+\epsilon,\tau_1]$. In particular for $\tau_1$}, which contradicts the definition of $\tau_1$ in \eqref{contrtau1} and proves the non strict inequality \eqref{nonstrct}.

Now we show the strict inequality \eqref{mumonot}. Suppose on the contrary that there exists a $\tau_4\in\mathbb{R}$ such that $\mu(\tau_4)=\tilde{\mu}(\tau_4)$. 

Choose $\tau^*<\tau_4$ as before, such that the strict inequality \eqref{lintaustar} holds for $\tau\in(-\infty,\tau^*]$. Due to the non-strict inequality \eqref{nonstrct}, we know $\mu(\tau)\geq\tilde{\mu}(\tau)$ for $\tau\in(\tau^*,\tau_4)$.   
Integrate backwards from $\tau_4$ to $\tau^*$, by reversing the orientation of $\tau\in[\tau^*,\tau_4]$ through $\tilde{\tau}:=-\tau$, so that $\tilde{\tau}\in[\tau_4,\tau^*]$.
 The integral formulation of the ODE yields
\begin{equation*}
   \mu(\tau^*)-\mu(\tau_3)= \int_{\tau_3}^{\tau^*} F_2(\lambda,\theta,\rho,\mu) d\tilde{\tau}
\end{equation*} 
with similar equation for $\tilde{\mu}$.

Hence, the same methods from equations \eqref{CUBA} and \eqref{CUBA2} can be applied for the difference ${\mu}({\tau^*})-\tilde{\mu}({\tau^*})$, yielding the inequality ${\mu}({\tau}^*)-\tilde{\mu}({\tau}^*)\leq 0$, which contradicts the definition of $\tau^*$. This proves the inequality \eqref{mumonot}.

Analogously, the above arguments can be used to prove the monotonicity in the radial coordinate. There are two mild adaptations in the proof. Firstly, one does not need to study the linearized flow for the radius, since the initial data is already ordered by $d<\tilde{d}$. Secondly, to obtain \eqref{rhoormu}, one needed to bound $|\tilde{\rho}-\rho|$. Here, we need to bound $|\tilde{\mu}-\mu|$, which is a continuous function on the compact interval $[\tau_2,\tau_1]$ and hence attains a maximum. Then, the mean value theorem is used for $|\tilde{\rho}-\rho|$.

The monotonicity in the parameter $\lambda$ is seen by comparing the flow \eqref{shootCI} as $\lambda$ increases.
\begin{flushright}
	$\blacksquare$
\end{flushright}
\end{pf}

\subsection{Intersection of shooting curves: finding equilibria}
The shooting curves $M^u_{\pi/2}$ and $M^s_{\pi/2}$ intersect at the constant equilibria $d\in \{-1,0,1\}$. 

If $d>1$, the shooting curves $M^u_{\pi/2}$ and $M^s_{\pi/2}$ are monotone in the initial data $d$ and $e$, respectively. Indeed, the former increase in the $p$ direction as $\theta$ increases, whereas the latter decreases in the $p$ direction, for any $\lambda\in \mathbb{R}_+$. Hence, they do not intersect. Analogously for $d<1$.

Consider the case that $|d|<1$. We show that intersections of the shooting curves only occur either at the $u$ or $p$-axis. Then we show how many intersections there are with those axis.
\begin{lem}
    $M^u \cap M^s \subseteq \{ (\theta,u,p)\in\mathbb{R}^3 \text{ $|$ }p=0 \text{ or } u=0\}$.
\end{lem}
\begin{pf}
Towards a contradiction, suppose there is an intersection point $(u,p)\in M^u \cap M^s$ which is not in these axis.

If $(u,p)\in M^u$, then $(-u,-p)\in M^u$, due to reflection symmetry \eqref{sym2}. Similarly, if $(u,p)\in M^s$, then $(-u,-p)\in M^s$. Therefore, 
\begin{equation*}
    (-u,-p)\in M^u \cap M^s.
\end{equation*}

Also, if $(u,p)\in M^u$, then $(u,-p)\in M^s$, due to the construction of $M^s$, which is done by the time reversal \eqref{sym1} of $M^u$. Similarly, if $(u,p)\in M^s$, then $(u,-p)\in M^u$, i.e.
\begin{equation*}
    (u,-p)\in M^u \cap M^s.
\end{equation*} 

The same arguments in the above two paragraph using both the time reversal \eqref{sym1} and reflection symmetry \eqref{sym2} yield
\begin{equation*}
    (-u,p)\in M^u \cap M^s.
\end{equation*}

Therefore there are four points with the same radius in the intersection $M^u \cap M^s$ and none of those lie in the $u$ or $p$-axis. The pigeon hole principle guarantees that at least two of those four points were constructed with initial data $d$ either in $(0,1)$ or $(-1,0)$, contradicting the monotonicity of the radius \eqref{rhomonot}, and proving the lemma.
\begin{flushright}
	$\blacksquare$
\end{flushright}
\end{pf}

The next step is to find exactly how many intersections there are between the stable and unstable shooting curves, for each $\lambda\in\mathbb{R}_+$. As $\lambda$ increases, the shooting curves change due to the continuous dependence on the parameter, yielding a different attractor. See \cite{Henry81} for the dependence of the attractor on parameters.

There are always three trivial equilibria $0,\pm1$ in the intersection of the shooting curves. A new pair of equilibria appears when $\lambda$ crosses an eigenvalue of the spherical laplacian $\lambda_k$. This characterizes the pitchfork bifurcations that occur at each $\lambda_k$ and gives a different proof of such results, as in \cite{Nascimento89}.
\begin{lem}\label{CHEGA}
    Consider $\lambda\in(\lambda_k,\lambda_{k+1})$, where $\lambda_k$ is the $k$-th eigenvalue of the axially symmetric Laplacian with $k\in\mathbb{N}_0$.
    
    Then there are $2k+3$ intersections of $M^u \cap M^s$, and the angle of the tangent vector of the unstable shooting curve at $(0,0)\in M^u_{\pi/2}$ is given by $\mu(\lambda_{k+1})=\frac{\pi}{2}(k+1)$.
\end{lem}
\begin{pf}
The proof follows by induction on $k\in\mathbb{N}_0$. For the basis of induction, $k=0$, it is proved that there are three equilibria for $\lambda\in (0,\lambda_1)$ and that $\mu(\lambda_1)=\frac{\pi}{2}$ at $(0,0)\in M^u_{\pi/2}$.

Consider the angle $\mu\left(\lambda,d\right)$ of the tangent vector of the unstable shooting curve $M^u_{\frac{\pi}{2}}$ with initial data $d\in [0,1]$ and $\lambda\in\mathbb{R}_+$.

For $\lambda_0=0$, the shooting flow \eqref{shootCI} implies that $p\equiv 0$ and hence the unstable shooting manifold is given by the $u$-axis. Therefore $\mu(0)=0$. By continuous dependence on $\lambda$, this curve changes a little for $\lambda$ small. Moreover, due to the monotonicities \eqref{taumonot}, \eqref{mumonot} and \eqref{rhomonot} for $d\in (0,1)$, the unstable shooting manifold spirals clockwise towards the trivial equilibria $0$. Considering the appropriate reflections through Symmetries 1 and 2, one obtains the full unstable and stable curves as below, for small $\lambda$.

We now show this shape persists as $\lambda$ is increased up to $\lambda_1$.
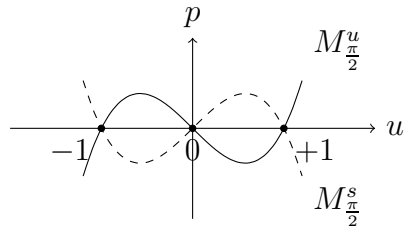
\begin{figure}[ht]\centering
\begin{tikzpicture}[scale=1.2]
    \draw[->] (-2,0) -- (2,0) node[right] {$u$};
    \draw[->] (0,-1) -- (0,1) node[above] {$p$};
    
    \filldraw [black] (-1,0) circle (1pt) node[anchor=north east]{$-1$};
    \filldraw [black] (1,0) circle (1pt) node[anchor=north west]{$+1$};
    \filldraw [black] (0,0) circle (1pt) node[anchor=north]{$0$};

    \draw [domain=-1.2:1.2,variable=\t,smooth] plot ({\t},{\t*(\t-1)*(\t+1)}) node[anchor=south west]{$M^u_{\frac{\pi}{2}}$};
    \draw [domain=-1.2:1.2,variable=\t,smooth, dashed] plot ({\t},{-\t*(\t-1)*(\t+1)}) node[anchor=north west]{$M^s_{\frac{\pi}{2}}$}; 
\end{tikzpicture}
\caption{Shooting curves of \eqref{shootCI} for $\lambda\in(0,\lambda_1)$}
\end{figure}

Recall that the angle $\mu$ is monotone in $\lambda$ for fixed $d\in \mathbb{R}$ as in \eqref{lambdamonot}. Moreover, $\lim_{\lambda\to\infty}\mu(\lambda,d)=\infty$, for any $\tau\in \mathbb{R}$ and $d\in(-1,1)$. Indeed, it follows by combining that \eqref{anglelambda} is increasing in $\lambda$, and \eqref{totangle} with the symmetry \eqref{sym1}, which implies that the stable angle is minus the unstable angle. 

Therefore, there is a $\lambda_*>0$ such that $\mu(\lambda_*,0)=\frac{\pi}{2}$. We have to prove that $\lambda_*=\lambda_1$ and that there are no new equilibria for $\lambda\in (0,\lambda_*)$.  

For fixed $\lambda$, the angle monotonicity \eqref{mumonot} implies that the biggest value that $\mu$ attains is at $d=0$. Together with the monotonicity in $\lambda$, we have that
\begin{equation}\label{enoughofthisthesis}
    \mu(\lambda,d)<\mu(\lambda_*,d)<\mu(\lambda_*,0)=\frac{\pi}{2} 
\end{equation} 
for $\lambda<\lambda_*$ and $d\in(0,1)$. Hence, there is no intersection of the unstable shooting curve with the negative $p$-axis, described by the angle $\pi/2$ in polar coordinates. 

Therefore, there is also no intersection of the unstable shooting curve with the negative $u$-axis (given by $\pi$ in polar coordinates). This occurs since the angle $\mu$ is continuous and monotone, hence the shooting curve would have to cross the negative $p$-axis at polar angle $\pi/2$, which was shown that does not occur in \eqref{enoughofthisthesis}. Similarly for the positive $p$-axis (described by the angle $3\pi/2$ in polar coordinates), and positive $u$-axis (at polar angle $0$ or $2\pi$).

The construction of the remaining part of the unstable manifold $M^u_{\pi/2}$ for ${d\in (-1,0)}$, through symmetry \eqref{sym2}, implies there is no intersection of this piece of the unstable shooting curve with the $u$ or $p$-axis. Hence, the only intersection points of the unstable shooting curve with the $u$ or $p$-axis lie in the trivial equilibria $d=-1,0,1$. 

Moreover, due to symmetry \eqref{sym1} and the construction of the shooting stable manifold $M^s_{\pi/2}$, there are no intersection points of the shooting stable manifold with the $u$ or $p$-axis, except at the equilibria with initial data $e=-1,0,1$. 

This proves there are no other equilibria for $\lambda\in (0,\lambda_*)$. We now show $\lambda_*=\lambda_1$. 

Due to the symmetry \eqref{sym1}, the angle of the tangent of the manifold $M^s_{\pi/2}$ at $e=0$ will be $-\mu(\lambda_*,0)=-\frac{\pi}{2}$. Hence, the angle between those tangent vectors is $\pi$, as in \eqref{totangle}. This occurs exactly when $\lambda_*=\lambda_1$, as the definition of the eigenvalue $\lambda_1$ through the eigenvalue problem in polar coordinates \eqref{EVpolar}. This proves the basis of induction.  

For the induction step, suppose that for $\lambda\in(\lambda_{k-1},\lambda_{k})$, there are $2(k-1)+3$ equilibria and $\mu(\lambda_{k})=\frac{\pi}{2}k$. Note the last condition informs how many times the unstable shooting curve has crossed the $u$ and $p$ axis. We shall prove that for $\lambda\in(\lambda_{k},\lambda_{k+1})$, there are $2k+3$ equilibria, $\mu(\lambda_{k+1})=\frac{\pi}{2}(k+1)$ and $\lambda_{k+1}$ is the $(k+1)$-th eigenvalue.

There exists a $\lambda^*>\lambda_k$ such that $\mu(\lambda^*,0)=\frac{\pi}{2}(k+1)$, due to the monotonicity in $\lambda$ as \eqref{lambdamonot}. The arguments to show that $\lambda^*=\lambda_{k+1}$ and that two new equilibria appear for $\lambda\in (\lambda_k,\lambda^*)$ are analogous as the basis of induction. 

There are two cases, depending on the parity of $k$. This influences which axis the shooting curve intersects and where the new equilibria appear, as $\lambda$ croses $\lambda_k$.

If $k$ is odd, then the new equilibria appear in the $p$-axis. We illustrate such case in the figure below, when $\lambda$ crosses $\lambda_1$.
\begin{figure}[ht]\centering
\begin{tikzpicture}[scale=0.4]
    \draw[->] (-6,0) -- (6,0) node[right] {$u$};
    \draw[->] (0,-3) -- (0,3) node[above] {$p$};

    \filldraw [black] (-3.14,0) circle (3pt) node[anchor=north east]{$-1$};
    \filldraw [black] (3.14,0) circle (3pt) node[anchor=north west]{$+1$};
    \filldraw [black] (0,0) circle (3pt) node[anchor=north]{$0$};
    \filldraw [black] (0,-1.57) circle (3pt);
    \filldraw [black] (0,1.57) circle (3pt);
    
    \draw [domain=0:3.14,variable=\t,smooth] plot ({\t*cos(\t r)},{\t*sin(\t r)}); 
    \draw [domain=-3.6:-3.14,variable=\t,smooth] plot ({\t},{-(0.7)*(\t-3.14)*(\t+3.14)}); 
    \draw [domain=0:3.14,variable=\t,smooth] plot ({-\t*cos(\t r)},{-\t*sin(\t r)}); 
    \draw [domain=3.14:3.6,variable=\t,smooth] plot ({\t},{(0.7)*(\t-3.14)*(\t+3.14)})node[anchor=south west]{$M^u_{\frac{\pi}{2}}$}; 

    \draw [domain=0:3.14,variable=\t,smooth,dashed] plot ({\t*cos(\t r)},{-\t*sin(\t r)});
    \draw [domain=-3.6:-3.14,variable=\t,smooth,dashed] plot ({\t},{(0.7)*(\t-3.14)*(\t+3.14)}); 
    \draw [domain=0:3.14,variable=\t,smooth,dashed] plot ({-\t*cos(\t r)},{\t*sin(\t r)}); 
    \draw [domain=3.14:3.6,variable=\t,smooth,dashed] plot ({\t},{-(0.7)*(\t-3.14)*(\t+3.14)})node[anchor=north west]{$M^s_{\frac{\pi}{2}}$};  
\end{tikzpicture}
\caption{Shooting curves of \eqref{shootCI} for $\lambda\in(\lambda_1,\lambda_2)$}
\end{figure}

This can be proved as follows. By the induction hypothesis, we know $\mu(\lambda_k)=\frac{\pi}{2}k$, which means that the tangent of the shooting $M^u_{\pi/2}$ at $d=0$ is parallel to the $p$-axis for odd $k$. Since $\mu$ is increasing in $\lambda$ and the shooting curve is continuous, then the shooting curve nearby $d=0$ moves from the quadrants $\{ p>0,u<0\}$ and $\{ p<0,u>0\}$ to its compliment, as $\lambda$ crosses $\lambda_k$. This creates two new intersections of the unstable shooting curve with the $p$-axis. Due to the construction of the stable shooting curve $M^s_{\pi/2}$, it also intersects the $p$-axis in the same points.

Then, one repeat the arguments in the induction step in order to show there is no intersection of the shooting curves with the $u$-axis.

The only remaining claim to be proven is that equilibria can't disappear, after they appear. The only possibility for this to happen is if two equilibria within the $u$ or $p$ axis collide. Note that neighboring equilibria come from different parts of the initial data: either $d$ is in $(0,1)$ or $(-1,0)$. Hence, if they collide, it contradicts uniqueness of the shooting flow \eqref{shootCI}, since their initial data is different.

Analogously, the case when $k$ is even yields new equilibria in the $u$-axis. We illustrate this in the example below, as $\lambda$ crosses $\lambda_2$ and two new equilibria appear in the $u$-axis.

\begin{figure}[ht]\centering
    \begin{tikzpicture}[scale=0.3]
    \draw[->] (-10,0) -- (10,0) node[right] {$u$};
    \draw[->] (0,-5) -- (0,5) node[above] {$p$};
    
    \filldraw [black] (-4.71,0) circle (4pt) node[anchor=north east]{$-1$};
    \filldraw [black] (4.71,0) circle (4pt) node[anchor=north west]{$+1$};
    \filldraw [black] (0,0) circle (4pt) node[anchor=north]{$0$};
    \filldraw [black] (0,-3.14) circle (4pt);
    \filldraw [black] (0,3.14) circle (4pt);
    \filldraw [black] (-1.57,0) circle (4pt);
    \filldraw [black] (1.57,0) circle (4pt);
    
    \draw [domain=0:4.71,variable=\t,smooth] plot ({\t*sin(\t r)},{-\t*cos(\t r)}); 
    \draw [domain=-5.2:-4.71,variable=\t,smooth] plot ({\t},{-(0.7)*(\t-4.71)*(\t+4.71)}); 
    \draw [domain=0:4.71,variable=\t,smooth] plot ({-\t*sin(\t r)},{\t*cos(\t r)}); 
    \draw [domain=4.71:5.2,variable=\t,smooth] plot ({\t},{(0.7)*(\t-4.71)*(\t+4.71)})node[anchor=south west]{$M^u_{\frac{\pi}{2}}$}; 

    \draw [domain=0:4.71,variable=\t,smooth,dashed] plot ({\t*sin(\t r)},{\t*cos(\t r)}); 
    \draw [domain=-5.2:-4.71,variable=\t,smooth,dashed] plot (-{\t},{-(0.7)*(\t-4.71)*(\t+4.71)}); 
    \draw [domain=0:4.71,variable=\t,smooth,dashed] plot ({-\t*sin(\t r)},{-\t*cos(\t r)}); 
    \draw [domain=4.71:5.2,variable=\t,smooth,dashed] plot (-{\t},{(0.7)*(\t-4.71)*(\t+4.71)}); 
    \end{tikzpicture}
    \caption{Shooting curves of \eqref{shootCI} for $\lambda\in(\lambda_2,\lambda_3)$}
\end{figure}
\begin{flushright}
	$\blacksquare$
\end{flushright}
\end{pf}
 

\subsection{Hyperbolicity: all intersections are transverse.} 

It is enough to check that the unstable and stable manifolds $M^u$ and $M^s$ are not tangent to the $u$ or $p$-axis. Indeed, $M^u$ is tangent to $u$-axis if, and only if $M^s$ is also, since one is obtained from the other through to the reflection $p\mapsto -p$. Similarly, $M^u$ is tangent to $p$-axis if, and only if there is another tangency of $M^u$ with the $p$-axis, due to the rotation $(u,p)\mapsto (-u,-p)$. Moreover, $M^s$ is obtained from $M^u$ through the reflection $p\mapsto -p$, hence $M^s$ is also tangent to the $p$-axis at the same points.
 
Recall that the tangent vector of the unstable shooting manifold is given by $(u_{d},p_{d})$ and satisfies the equation \eqref{tangentshoot}. Hence, this vector is tangent to the $p$-axis if it is vertical, that is, if the coordinate $u_{d}=0$. On the other hand, the coordinate in polar coordinates is $u=\rho \cos(\mu)$, and the chain rule implies 
\begin{equation*}
    0=u_{d}=(\rho \cos(\mu))_{d}=\rho_{d} \cos(\mu)-\rho \sin(\mu)\mu_{d}.
\end{equation*}

Algebraic manipulation yields $\mu= \arctan( \rho_{d} / (\rho \mu_{d}))$. Note that the monotonicity properties \eqref{mumonot} and \eqref{rhomonot} implies $\rho_{d}$ and $\mu_{d}$ are nonzero with different signs, for both cases that $d$ is either in $ (0,1)$ or $(-1,0)$. Moreover, the radius $\rho>0$. Therefore, the argument $\rho_{d}/(\rho\mu_{d})$ is strictly negative, and hence $\mu\in (-\pi/2,0)$. That is, the point where the tangency occurs is neither at the $u$, nor the $p$-axis, because those in polar coordinates are given by $\mu=\frac{\pi}{2}k$. This contradicts that intersections must occur at the $p$-axis.

Similarly a tangency occurs at the $u$-axis, if the vector is horizontal, namely $p_{d}=0$. In polar coordinates $p=-\rho \sin (\mu)$, the tangency condition and the chain rule implies that $\mu=\arctan(-\rho_{d}/(\rho\mu_{d}))$. By a similar analysis as above, the argument is strictly positive and hence $\mu\in (0,\pi/2)$. That is, the intersection does not occur in the $u$-axis, yielding a contradiction.

Therefore, all nontrivial equilibria are hyperbolic. The only equilibrium that can be non-hyperbolic is the trivial one when $d=0$, because in this case $\rho=0$ and hence the tangent is parallel to the $u$ or $p$-axis. Indeed, this was proved to be the case in Lemma \ref{CHEGA}, whenever $\lambda=\lambda_k$ for all $k\in\mathbb{N}_0$.

\subsection{Obtaining the permutation.}

We construct the Sturm permutation for $\lambda \in (\lambda_{k},\lambda_{k+1})$ by induction on $k\in\mathbb{N}_0$. The idea is to label the intersections of the unstable and stable manifolds, firstly along $M^u$ following its parametrization given by the initial data $d$ from $-\infty$ to $\infty$, and secondly along $M^s$ following its parametrization given by $e$ from $-\infty$ to $\infty$. 

For $k=0$, that is, $\lambda \in(\lambda_0,\lambda_1)$, there are no other intersections of the shooting curves, except the trivial equilibria $d=-1,0,1$. Noticing how the shooting curve was constructed before, this is exactly their order along both $M^u$ and $M^s$ following their parametrization $d,e$ from $-\infty$ to $\infty$. Hence, the permutation is the identity $\sigma=id$, since their order is the same along both $M^u$ and $M^s$.

For the induction step, we find the permutation for $\lambda \in (\lambda_{k},\lambda_{k+1})$, with $k\geq 1$, supposing that the permutation for $\lambda \in (\lambda_{k-1},\lambda_{k})$ is given by 
\begin{equation}\label{perM}
    \sigma= (2, 2k) (4, 2k-2) ... 
\end{equation}
where $(j, l)$ is a transposition in the group of permutations $S_N$ with appropriate $N$. 

For $k\geq 1$ and $\lambda<\lambda_{k}$ with small $|\lambda-\lambda_{k}|$, there are 
$N=2k+1$ equilibria as in Lemma \ref{CHEGA}, and hence $\lceil k/2\rceil$ transpositions in the Sturm permutation, where $\lceil . \rceil$ denotes the ceiling function. Note $\sigma$ contains all even numbers less or equal $N=2k+1$.

There are two cases: either $k$ is even or odd. The previous construction of the shooting curve implies that it rotates clock-wise around the trivial equilibria $0$. The parity of $k$ tells how the shooting curve behaves for $\lambda<\lambda_{k}$, in particular, if the equilibria nearby the trivial equilibria $0$ is obtained by an intersection with the $u$ or $p$ axis.

Suppose $k$ is odd. Labeling the equilibria along $M^u$ and $M^s$ for $\lambda<\lambda_k$, the trivial equilibria $0$ is labeled $k+1$, since there are $k$ equilibria before it along the unstable manifold. Hence, the nearby equilibria are labeled by $k$ for the equilibria before, and $k+2$ for the equilibria after it. Moreover, the last transposition in the permutation \eqref{perM} is $(k+1,k+1)$, since $k+1$ is even. 

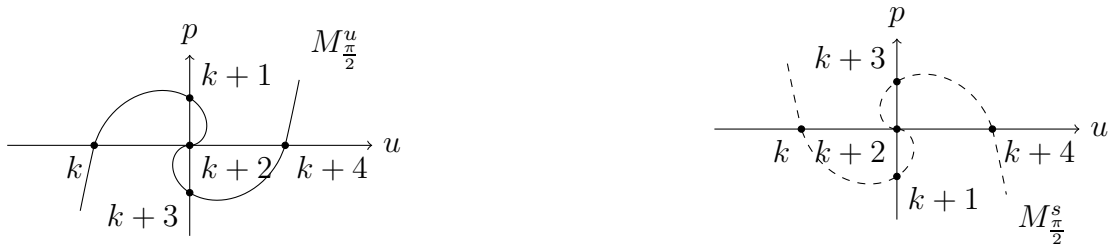
\begin{figure}[!htb]\label{fig:bef}\centering
    \begin{tikzpicture}[scale=1.2]
    \draw[->] (-2,0) -- (2,0) node[right] {$u$};
    \draw[->] (0,-1) -- (0,1) node[above] {$p$};
    
    \filldraw [black] (-1,0) circle (1pt) node[anchor=north east]{$k$};
    \filldraw [black] (1,0) circle (1pt) node[anchor=north west]{$k+2$};
    \filldraw [black] (0,0) circle (1pt) node[anchor=north east]{$k+1$};

    \draw [domain=-1.2:1.2,variable=\t,smooth] plot ({\t},{\t*(\t-1)*(\t+1)}) node[anchor=south west]{$M^u_{\frac{\pi}{2}}$};
    \draw [domain=-1.2:1.2,variable=\t,smooth, dashed] plot ({\t},{-\t*(\t-1)*(\t+1)}) node[anchor=north west]{$M^s_{\frac{\pi}{2}}$}; 
\end{tikzpicture}
\caption{Labeling of equilibria for $\lambda<\lambda_k$ with $k$ odd}
\end{figure}

The labeling within the unstable manifold for the equilibria labeled less than $k+1$ will not change. Moreover, as $\lambda$ cross $\lambda_k$, two new equilibria appear, one on each side along the unstable manifold. The trivial equilibria 0, which was labeled $k+1$ for $\lambda<\lambda_k$, will be shifted by $1$, yielding $k+2$ for $\lambda>\lambda_k$. All other having label bigger than $k+1$ will be shifted by two. See the figure below.

Similarly occur a change of labeling along the stable manifold.

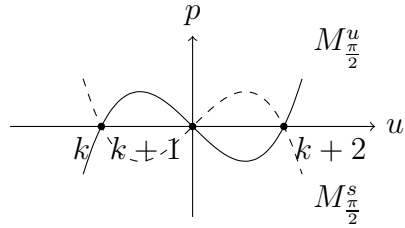
\begin{figure}[!htb]
\minipage{0.4\textwidth}
    \begin{tikzpicture}[scale=0.4]
    \draw[->] (-6,0) -- (6,0) node[right] {$u$};
    \draw[->] (0,-3) -- (0,3) node[above] {$p$};

    \filldraw [black] (-3.14,0) circle (3pt) node[anchor=north east]{$k$};
    \filldraw [black] (3.14,0) circle (3pt) node[anchor=north west]{$k+4$};
    \filldraw [black] (0,0) circle (3pt) node[anchor=north west]{$k+2$};
    \filldraw [black] (0,-1.57) circle (3pt)node[anchor=north east]{$k+3$};
    \filldraw [black] (0,1.57) circle (3pt)node[anchor=south west]{$k+1$};
    
    \draw [domain=0:3.14,variable=\t,smooth] plot ({\t*cos(\t r)},{\t*sin(\t r)}); 
    \draw [domain=-3.6:-3.14,variable=\t,smooth] plot ({\t},{-(0.7)*(\t-3.14)*(\t+3.14)}); 
    \draw [domain=0:3.14,variable=\t,smooth] plot ({-\t*cos(\t r)},{-\t*sin(\t r)}); 
    \draw [domain=3.14:3.6,variable=\t,smooth] plot ({\t},{(0.7)*(\t-3.14)*(\t+3.14)})node[anchor=south west]{$M^u_{\frac{\pi}{2}}$}; 
    \end{tikzpicture}
\endminipage\hfill
\minipage{0.4\textwidth}%
    \begin{tikzpicture}[scale=0.4]
    \draw[->] (-6,0) -- (6,0) node[right] {$u$};
    \draw[->] (0,-3) -- (0,3) node[above] {$p$};

    \filldraw [black] (-3.14,0) circle (3pt) node[anchor=north east]{$k$};
    \filldraw [black] (3.14,0) circle (3pt) node[anchor=north west]{$k+4$};
    \filldraw [black] (0,0) circle (3pt) node[anchor=north east]{$k+2$};
    \filldraw [black] (0,-1.57) circle (3pt)node[anchor=north west]{$k+1$};
    \filldraw [black] (0,1.57) circle (3pt)node[anchor=south east]{$k+3$};
    
    \draw [domain=0:3.14,variable=\t,smooth,dashed] plot ({\t*cos(\t r)},{-\t*sin(\t r)});
    \draw [domain=-3.6:-3.14,variable=\t,smooth,dashed] plot ({\t},{(0.7)*(\t-3.14)*(\t+3.14)}); 
    \draw [domain=0:3.14,variable=\t,smooth,dashed] plot ({-\t*cos(\t r)},{\t*sin(\t r)}); 
    \draw [domain=3.14:3.6,variable=\t,smooth,dashed] plot ({\t},{-(0.7)*(\t-3.14)*(\t+3.14)})node[anchor=north west]{$M^s_{\frac{\pi}{2}}$};  
\end{tikzpicture}
\endminipage
\caption{Labeling of equilibria for $\lambda>\lambda_k$ with $k$ odd.}
\end{figure}

One only has to check what happens to the permutation $\sigma$ as $\lambda$ crosses $\lambda_k$: a new transposition $(k+1, k+3)$ is added in the permutation. Note $k$ is odd, and hence both $k+1$ and $k+3$ are even. Moreover, the transposition, which was $(k+1,k+1)$ must be shifted to $(k+2,k+2)$, yielding the identity transposition and not changing \eqref{perM}. 

Therefore, the number of transpositions $\lceil k/2\rceil$ does not change. The only difference is the relabeling of equilibria within the permutation, described above, yielding the desired permutation. The case when the orientation of the unstable manifold is reversed have identical arguments.

For $k$ even, the above argument can be adapted. Notice that there are $2k+1$ equilibria, and again the trivial equilibria $0$ is labeled by $k+1$. As $\lambda$ crosses $\lambda_k$, there are two new equilibria along the unstable manifold. Hence, the ones before $k$ should not be relabeled, the origin $k+1$ for $\lambda<\lambda_k$ should be relabeled by $k+2$ for $\lambda>\lambda_k$, and all equilibria with label bigger than $k+1$ should be shifted by $2$.

Again, since the stable manifold is obtained by the reflection of the unstable manifold with respect to the $u$-axis, then one can see that the new permutation that should be added is $(k, k+4)$. Notice those are even numbers. Again, the transposition of the origin does not change \eqref{perM}, since it yields the identity transposition given by $(k+2,k+2)$. Similarly when the orientation of the unstable manifold is reversed.

\subsection{Obtaining the attractor.}

The permutation obtained above is the same as the regular Chafee-Infante problem. Hence, the attractors are connection-wise the same, since the conditions for the existence of heteroclinics are the same.

\section{Discussion}

The shooting method used to construct the attractor generalizes the bifurcation result in \cite{Nascimento89} for radially symmetric solutions in the disk. Indeed, not only we are able to prove the existence of bifurcating equilibria, but can also compute secondary bifurcations that might occur, hyperbolicity of all equilibria, their Morse indices and how they fit together in the attractor, by computing heteroclinic trajectories. 

After the construction of the Sturm attractor for the parabolic equation with singular coefficients, we see that if the Sturm permutation for the singular case coincide with the permutation for the case of regular coefficients, then the attractors for both cases coincide. This happens since they are both constructed in the same way, yielding the same necessary and sufficient conditions for heteroclinics to exist, regarding the zero numbers and Morse indices.

The construction of the Sturm attractor in the case of general singular diffusion is not proved here, but the above arguments can be replicated without severe modifications. 

\end{document}